\documentclass[11pt]{article}
\usepackage{latexsym,amssymb}
\usepackage{color}
\usepackage{amsmath}

\newtheorem{Theorem}{\bf Theorem}[section]
\newtheorem{Lemma}{\bf Lemma}[section]
\newtheorem{Proposition}{\bf Proposition}[section]
\newtheorem{Corollary}{\bf Corollary}[section]
\newtheorem{Remark}{\bf Remark}[section]
\newtheorem{Example}{\bf Example}[section]
\newtheorem{Definition}{\bf Definition}[section]

\newenvironment{theorem}{\begin{Theorem}$\!\!\!$}{\end{Theorem}}
\newenvironment{lemma}{\begin{Lemma}$\!\!\!$}{\end{Lemma}}

\newenvironment{remark}{\begin{Remark}$\!\!\!$}{\end{Remark}}

\numberwithin{equation}{section}

\begin{document}

\title{Asymptotic expansions of solutions\\ of fractional diffusion equations}
%\title{Asymptotics for fractional diffusion equations}
%
%
\author{Kazuhiro Ishige, 
Tatsuki Kawakami and 
Hironori Michihisa}
\date{}
\maketitle
%%%%%%%%%%%%%%%%%%%%%%%%%%%%%%%%%%%%%%
%%%%%%%%%%%%%%%%%%%%%%%%%%%%%%%%%%%%%%
\begin{abstract}
In this paper we obtain the precise description of the asymptotic behavior 
of the solution $u$ of 
$$
\partial_t u+(-\Delta)^{\frac{\theta}{2}}u=0\quad\mbox{in}\quad{\bf R}^N\times(0,\infty),
\qquad
u(x,0)=\varphi(x)\quad\mbox{in}\quad{\bf R}^N,
$$
where $0<\theta<2$ and $\varphi\in L_K:=L^1({\bf R}^N,\,(1+|x|)^K\,dx)$ with $K\ge 0$. 
Furthermore, we develop the arguments in \cite{IK02} and \cite{IKK02} 
and establish a method to obtain the asymptotic expansions 
of the solutions to a nonlinear fractional diffusion equation 
$$
\partial_t u+(-\Delta)^{\frac{\theta}{2}}u=|u|^{p-1}u\quad\mbox{in}\quad{\bf R}^N\times(0,\infty),
$$
where $0<\theta<2$ and $p>1+\theta/N$.  
\end{abstract}
\vspace{10pt}

\vspace{40pt}
\noindent Addresses:

\smallskip
\noindent K. I.:  Mathematical Institute, Tohoku University,
Aoba, Sendai 980-8578, Japan.\\
\noindent 
E-mail: {\tt ishige@math.tohoku.ac.jp}\\

\smallskip
\noindent 
T. K.: Department of Mathematical Sciences, Osaka Prefecture University, 
Sakai 599-8531, Japan. \\
\noindent 
E-mail: {\tt kawakami@ms.osakafu-u.ac.jp}\\

\smallskip
\noindent 
H. M.: Mathematical Institute, Tohoku University,
Aoba, Sendai 980-8578, Japan.\\
\noindent 
E-mail: {\tt hironori.michihisa.p1@dc.tohoku.ac.jp}\\
\newpage
%%%%%%%%%%%%%%%%%%%%%%%%%%%%%%%%%%%%%
%%%%%%%%%%%%%%%%%%%%%%%%%%%%%%%%%%%%%
\section{Introduction}
%%%%%%%%%%%%%%%%%%%%%%%%%%%%%%%%%%%%%
%%%%%%%%%%%%%%%%%%%%%%%%%%%%%%%%%%%%%
This paper is concerned with the Cauchy problem 
for a nonlinear fractional diffusion equation 
\begin{equation}
\label{eq:1.1}
\partial_t u+(-\Delta)^{\frac{\theta}{2}}u=|u|^{p-1}u\quad\mbox{in}\quad{\bf R}^N\times(0,\infty),
\quad
u(x,0)=\varphi(x)\quad\mbox{in}\quad{\bf R}^N,
\end{equation}
where $N\ge 1$, $\partial_t:=\partial/\partial t$, 
$(-\Delta)^{\theta/2}$ is the fractional power of the Laplace operator with $0<\theta<2$, 
$p>1+\theta/N$ and $\varphi\in L^\infty({\bf R}^N)\cap L^1({\bf R}^N)$.  
We say that 
a continuous function $u$ in ${\bf R}^N\times(0,\infty)$ is a solution of \eqref{eq:1.1} 
if $u$ satisfies 
\begin{equation}
\label{eq:1.2}
u(x,t)=\int_{{\bf R}^N}G_\theta(x-y,t)\varphi(y)\,dy
+\int_0^t\int_{{\bf R}^N}G_\theta(x-y,t-s)F(u(y,s))\,dy\,ds
\end{equation}
for $x\in{\bf R}^N$ and $t>0$, where $F(u):=|u|^{p-1}u$. 
Here $G_\theta=G_\theta(x,t)$ is the fundamental solution of 
the linear fractional diffusion equation
\begin{equation}
\label{eq:1.3}
\partial_t u+(-\Delta)^{\frac{\theta}{2}}u=0\quad\mbox{in}\quad{\bf R}^N\times(0,\infty). 
\end{equation} 
Problem~\eqref{eq:1.1} appears in the study of  nonlinear problems with anomalous diffusion 
and the Laplace equation with a dynamical boundary condition 
and it has been studied extensively by many mathematicians 
(see \cite{AF}, \cite{FIK}, \cite{FK}, \cite{HKN}, \cite{HI}, \cite{IK01}, \cite{IKK01}, \cite{IKK02}, \cite{S} and references therein). 
Among others, 
Sugitani~\cite{S} showed that, 
if $1<p\le 1+\theta/N$, then problem~\eqref{eq:1.1} possesses no positive global in time solutions. 
On the other hand, if $p>1+\theta/N$ and $\|\varphi\|_{L^{N(p-1)/\theta,\infty}({\bf R}^N)}$ is sufficiently small, 
then there exists a global in time solution $u$ of \eqref{eq:1.1} such that 
\begin{equation}
\label{eq:1.4}
\lim_{t\to\infty}t^{\frac{N}{\theta}(1-\frac{1}{q})}\|u(t)-m_uG_\theta(t)\|_{L^q({\bf R}^N)}=0
\end{equation}
for any $1\le q\le\infty$, where 
$$
m_u:=\lim_{t\to\infty}\int_{{\bf R}^N} u(x,t)\,dx. 
$$
(See \cite{IKK01} and \cite{IKK02}.)
In this paper we investigate the higher order asymptotic expansions of the solutions to \eqref{eq:1.1} 
satisfying \eqref{eq:1.4}. 
\vspace{5pt}

Let us consider the Cauchy problem for a nonlinear heat equation 
\begin{equation}
\label{eq:1.5}
\partial_t v-\Delta v=f(x,t,v,\nabla v)\quad\mbox{in}\quad{\bf R}^N\times(0,\infty),
\quad
v(x,0)=\varphi(x)\quad\mbox{in}\quad{\bf R}^N, 
\end{equation}
where $f$ is a continuous function in ${\bf R}^N\times(0,\infty)\times{\bf R}^{N+1}$ and 
$$
\varphi\in L_K:=L^1({\bf R}^N,(1+|x|)^K\,dx)
\quad\mbox{for some $K\ge 0$.}
$$
If a classical solution $v$ of \eqref{eq:1.5} satisfies 
$$
\left|f(x,t,v(x,t),(\nabla v)(x,t))\right|\le C(1+t)^{-A}\left[|v(x,t)|+(1+t)^{\frac{1}{2}}|(\nabla v)(x,t)|\right]
$$
in ${\bf R}^N\times(0,\infty)$
for some $C>0$ and $A>1$, then 
the solution~$v$ behaves like a suitable multiple of the Gauss kernel as $t\to\infty$. 
More precisely,   
\begin{equation}
\label{eq:1.6}
\lim_{t\to\infty}\,t^{\frac{N}{2}(1-\frac{1}{q})}\|v(t)-m_v\Gamma(t)\|_{L^q({\bf R}^N)}=0
\end{equation}
holds for any $1\le q\le\infty$, 
where 
$$
\Gamma(x,t):=(4\pi t)^{-\frac{N}{2}}\exp\left(-\frac{|x|^2}{4t}\right),
\qquad
m_v:=\lim_{t\to\infty}\int_{{\bf R}^N}v(x,t)\,dx. 
$$
See e.g., \cite[Theorem~3.1]{IK03}. 
Generally, the higher order asymptotic expansions of the solutions 
depend on the decay of the solutions at the space infinity 
and they have been studied in many papers, 
see e.g., 
\cite{DK}, \cite{FM}, \cite{HKN}, \cite{IIK}, \cite{IK01}, \cite{IK02}, \cite{IK03}, \cite{KP} and \cite{T}. 
Among others, 
the first and the second authors of this paper developed the arguments in \cite{IIK} 
and established a systematic method to obtain the higher order asymptotic expansions 
of the solutions satisfying \eqref{eq:1.6}
(see \cite{IK02} and \cite{IK03}). 
Their arguments are based on the fact that $v(t)\in L_K$ for any $t>0$. 

We consider fractional diffusion equations. 
Let $0<\theta<2$ and set
$$
\big[S_\theta(t)\varphi\big](x)
:=\int_{{\bf R}^N}G_\theta(x-y,t)\varphi(y)\,dy,\qquad x\in{\bf R}^N,\,t>0. 
$$
Since $G_\theta(\cdot,t)\in L_K$ for $t>0$ if and only if $K<\theta$, 
we see that 
$$
\left\{S_\theta(t)\varphi\,:\,\varphi\in L_K\right\}\subset L_K\quad\mbox{for any $t>0$}
\quad\mbox{if and only if}\quad
K<\theta.
$$
This means that, in the case $K\ge\theta$, 
$L_K$ is not a suitable function space of initial functions 
for the fractional heat semigroup $\{S_\theta(t)\}_{t\ge 0}$. 
As far as we know, 
there are no systematic methods 
which are available to the study of the asymptotic expansions of the solutions to \eqref{eq:1.1} 
with $\varphi\in L_K$ in the case $K\ge\theta$. 
For the case $0\le K<\theta$, see \cite{IKK02}. 
(For related results, see e.g., \cite{Iwa}, \cite{OY}, \cite{Y03}, \cite{YS} and references therein.) 
\vspace{5pt}

In this paper we refine and develop the arguments in \cite{IKK02} and 
establish a systematic method to obtain the asymptotic expansions of the solutions to \eqref{eq:1.1} 
without any restrictions such as $K<\theta$. 
More precisely, the purposes of this paper are the following: \\

Let $0<\theta<2$, $j\in\{0,1,2,\dots\}$ and $K\ge 0$.\begin{itemize}
 \item[{\rm (i)}] 
We obtain the precise description of the asymptotic behavior of $\nabla^j [S_\theta(t)\varphi]$ with $\varphi\in L_K$ 
and identify the largest subspace $X$ of $L^1({\bf R}^N)$ satisfying 
\begin{itemize}
 \item
 $X\subset L_K$,
 \item 
 $\left\{\nabla^j[S_\theta(t)\varphi]\,:\,\varphi\in X\right\}\subset X$ for any $t>0$.
\end{itemize}
  See Theorem~\ref{Theorem:1.1} and Theorem~\ref{Theorem:1.2}~(a);
  \item[{\rm (ii)}] 
  We find a subspace $Y$ of $L^1({\bf R}^N)$ such that 
  \begin{itemize}
  \item 
  $L_K\subset Y$ and 
  the dimension of the quotient space $Y/L_K$ is finite,
  \item 
  $\displaystyle{\left\{\nabla^j [S_\theta(t)\varphi]\,:\,\varphi\in Y\right\}\subset Y}$ for any $t>0$.
  \end{itemize}
  See Theorem~\ref{Theorem:1.2}~(b);
  \item[{\rm (iii)}] 
  Let $\psi$ be a radially symmetric smooth function in ${\bf R}^N$ such that 
  $\partial_x^\alpha\psi\in L^\infty({\bf R}^N)\cap L_{K+|\alpha|}$ for $\alpha\in {\bf M}_K$. 
  By using $S_\theta(t)\psi$ and its derivatives, 
  we obtain the precise description of the asymptotic behavior of the solution~$u$ of \eqref{eq:1.1} 
  satisfying \eqref{eq:1.4}.
  See Sections~3 and 4. 
\end{itemize}

We introduce some notation. 
For any $k\ge 0$, let $[k]\in{\bf N}\cup\{0\}$ be such that $k-1<[k]\le k$.
For any multi-index $\alpha\in{\bf M}:=({\bf N}\cup\{0\})^N$, 
set
$$
|\alpha|:=\displaystyle{\sum_{i=1}^N}\alpha_i,\quad
\alpha!:=\prod_{i=1}^N\alpha_i!,\quad
x^\alpha:=\prod_{i=1}^N x_i^{\alpha_i},\quad
\partial_x^\alpha:=
\frac{\partial^{|\alpha|}}{\partial x_1^{\alpha_1}\cdots\partial x_N^{\alpha_N}}.
$$
Let ${\bf M}_k:=\left\{\alpha\in{\bf M}:|\alpha|\le k\right\}$ if $k\ge 0$ and ${\bf M}_k:=\emptyset$ if $k<0$. 
For any $\alpha=(\alpha_1,\dots,\alpha_N)$, $\beta=(\beta_1,\dots,\beta_N)\in{\bf M}$, 
we say 
$\alpha\le\beta$
if $\alpha_i\le\beta_i$ for all $i\in\{1,\dots,N\}$. 
Let $\nabla:=(\partial/\partial x_1,\dots,\partial/\partial x_N)$. 
For any $\alpha\in{\bf M}$ and $0<\theta<2$, we define
\begin{equation}
\label{eq:1.7}
g_{\alpha,\theta}(x,t):=\frac{(-1)^{|\alpha|}}{\alpha!}(\partial_x^\alpha G_\theta)(x,t)
=\frac{(-1)^{|\alpha|}}{\alpha!}t^{-\frac{N+|\alpha|}{\theta}}(\partial_x^\alpha G_\theta)(t^{-\frac{1}{\theta}}x,1).
\end{equation}
(See ({\bf G})-(i) in Section~2.) 
For any $1\le q\le\infty$, let $\|\cdot\|_q$ be the usual norm of $L^q:=L^q({\bf R}^N)$. 
Furthermore, for any $\ell\ge 0$, 
we set 
$$
|||f|||_\ell:=\int_{{\bf R}^N}|x|^\ell|f(x)|\,dx,\qquad f\in L_\ell. 
$$
For a real vector space $V$ and a finite set $\{f_i\}_{i=1}^m$ in $V$,  
let
$$
\sum_{i=1}^m\,\langle f_i\rangle:=\biggr\{\sum_{i=1}^m c_if_i\,:\,\{c_i\}_{i=1}^m\subset{\bf R}\biggr\}\subset V.
$$
\vspace{3pt}

Now we are ready to state the main results of this paper. 
In Theorem~\ref{Theorem:1.1} we obtain  
the precise description of the asymptotic behavior of $S_\theta(t)\varphi$ with $\varphi\in L_K$. 
\begin{theorem}
\label{Theorem:1.1} 
Let $N\ge 1$, $0<\theta<2$, $j\in\{0,1,2,\dots\}$ and $K\ge 0$. 
For any $\varphi\in L_K$, set 
\begin{equation}
\label{eq:1.8}
v(x,t):=\big[S_\theta(t)\varphi\big](x)-\sum_{|\alpha|\le K}
\biggr[\int_{{\bf R}^N}x^\alpha \varphi(x)\,dx\biggr]g_{\alpha,\theta}(x,t). 
\end{equation}
Then 
\begin{equation}
\label{eq:1.9}
(\nabla^j v)(t)\in\left\{f\in L_K\,:\,\int_{{\bf R}^N}x^\alpha f(x)\,dx=0
 \,\,\,\mbox{for all $\alpha\in{\bf M}_K$}\right\}
\end{equation}
for $t>0$. 
Furthermore, the following hold:
\begin{itemize}
  \item[{\rm (a)}] 
  There exists a constant $C$ such that
  \begin{equation}
  \label{eq:1.10}
  t^{\frac{N}{\theta}(1-\frac{1}{q})}\|(\nabla^j v)(t)\|_q+
  t^{-\frac{\ell}{\theta}}|||(\nabla^j v)(t)|||_\ell\le Ct^{-\frac{K+j}{\theta}}|||\varphi|||_K,
  \quad t>0,
  \end{equation}
  for $1\le q\le\infty$ and $0\le\ell\le K$. 
  Here $C$ depends only on $N$, $\theta$, $j$ and $K$;
  \item[{\rm (b)}] For any $1\le q\le\infty$ and $0\le\ell\le K$,   
  $$
  \lim_{t\to\infty}t^{\frac{K+j}{\theta}}
  \left[t^{\frac{N}{\theta}(1-\frac{1}{q})}\|(\nabla^j v)(t)\|_q+
  t^{-\frac{\ell}{\theta}}|||(\nabla^j v)(t)|||_\ell\right]=0.
  $$
\end{itemize}
\end{theorem}
\begin{remark}
\label{Remark:1.1} 
{\rm (i)} Theorem~{\rm\ref{Theorem:1.1}} in the case $0\le K<\theta$ has been already proved in {\rm\cite{IKK02}}. 
However, the proof in this paper is simpler than that of {\rm\cite{IKK02}}.
\newline
{\rm (ii)} 
In the case $K\ge\theta$, 
some weaker estimates than \eqref{eq:1.10} were stated in {\rm\cite{Y03}} and {\rm\cite{YS}}. 
The proofs in {\rm\cite{Y03}} and {\rm\cite{YS}} 
were based on uncertain pointwise estimates 
of $(\partial_t^m\partial_x^\alpha G_\theta)(x,t)$,
where $m\in\{1,2,\dots\}$ and $\alpha\in{\bf M}$. 
\end{remark}
Next we characterize 
the behavior of $S_\theta(t)\varphi$ by the use of subspaces of $L^1$. 
\begin{theorem}
\label{Theorem:1.2} 
Let $N\ge 1$, $0<\theta<2$, $j\in\{0,1,2,\dots\}$ and $K\ge 0$. 
\begin{itemize}
  \item[{\rm (a)}] 
  The real vector space 
  $$
  {\mathcal L}_{K,j}^0:=\biggr\{f\in L_K\,:\,
  \int_{{\bf R}^N}x^\alpha f(x)\,dx=0\quad\mbox{for}\quad\alpha\in{\bf M}_{K-j-\theta}\biggr\}
  $$
  is the largest function space satisfying
  \begin{equation}
  \label{eq:1.11}
  {\mathcal L}_{K,j}^0\subset L_K
  \qquad\mbox{and}\qquad
  \bigcup_{t>0}\,\left\{\nabla^j[S_\theta(t)\varphi]\,:\,\varphi\in{\mathcal L}_{K,j}^0\right\}\subset{\mathcal L}_{K,j}^0.
  \end{equation} 
  \item[{\rm (b)}] 
  Define 
  $$
  {\mathcal L}_{K,j}:=L_K\,+\sum_{|\alpha|+m\theta\le K-j}\,
  \langle(\partial_t^m\partial_x^\alpha G_\theta)(\cdot,1)\rangle.
  $$
  Then
  $$
  \bigcup_{t>0}\,\left\{\nabla^j [S_\theta(t)\varphi]\,:\,\varphi\in{\mathcal L}_{K,j}\right\}
  \subset{\mathcal L}_{K,j}\subset L_{K'},
  $$
  where $0\le K'\le K$ with $K'<\theta$. 
\end{itemize}
\end{theorem}
Assertion~(b) implies that 
the dimension of the quotient space ${\mathcal L}_{K,j}/L_K$ is finite. 
Furthermore, as an application of Theorem~\ref{Theorem:1.1}, 
we study the behavior of the hot spots of $S_\theta(t)\varphi$. 
\begin{theorem}
\label{Theorem:1.3}
Let $N\ge 1$, $0<\theta<2$ and $K\ge 2$. 
Let $\varphi\in L_K$ and assume that 
$$
M(\varphi):=\int_{{\bf R}^N}\varphi(x)\,dx>0. 
$$
Then the hot spots of $S_\theta(t)\varphi$
$$
H(t):=\left\{x\in{\bf R}^N\,:\, [S_\theta(t)\varphi](x)=\sup_{y\in{\bf R}^N}[S_\theta(t)\varphi](y)\right\}
$$
is not empty for any $t>0$ and the following hold:
\begin{itemize}
  \item $H(t)$ consists of only one point $\{x(t)\}$ in finite time;
  \item $x(t)$ moves along a smooth curve in ${\bf R}^N$ and 
  $$
  \lim_{t\to\infty}x(t)=C(\varphi):=\int_{{\bf R}^N}x\varphi(x)\,dx\biggr/ M(\varphi). 
  $$
\end{itemize}
\end{theorem}
The movement of the hot spots is of independent interest 
and it has been studied for the heat equation, 
see \cite{CK}, \cite{FI}, \cite{IKabe}, \cite{IKoba} and references therein.
\vspace{3pt}
\par
The rest of this paper is organized as follows. 
In Section~2 we prove Theorem~\ref{Theorem:1.1} 
by refining the arguments in \cite{IIK} and \cite{IKK02}. 
Theorems~\ref{Theorem:1.2} and \ref{Theorem:1.3} follow from Theorem~\ref{Theorem:1.1}. 
In Section~3 we obtain the asymptotic expansions of the solutions to 
the fractional diffusion equation with an inhomogeneous term 
$$
\partial_t u+(-\Delta)^{\frac{\theta}{2}}u=F(x,t)\quad\mbox{in}\quad{\bf R}^N\times(0,\infty),
\qquad
u(x,0)=\varphi(x)\quad\mbox{in}\quad{\bf R}^N,
$$
where $\varphi\in L^\infty\cap L_K$ for some $K\ge 0$ and $F\in L^\infty(0,T:L_K)$ for any $T>0$. 
In Section~4 we apply the arguments in Sections~2 and 3 to obtain 
the asymptotic expansions of the solution of \eqref{eq:1.1} satisfying \eqref{eq:1.4}. 
Furthermore, 
combining the arguments in \cite{IK02}, 
we obtain the higher order asymptotic expansions of the solution~$u$ of \eqref{eq:1.1}. 
%%%%%%%%%%%%%%%%%%%%%%%%%%%%%%%%%%%%%
%%%%%%%%%%%%%%%%%%%%%%%%%%%%%%%%%%%%%
\section{Proofs of Theorems~\ref{Theorem:1.1}, \ref{Theorem:1.2} and \ref{Theorem:1.3}}
%%%%%%%%%%%%%%%%%%%%%%%%%%%%%%%%%%%%%
%%%%%%%%%%%%%%%%%%%%%%%%%%%%%%%%%%%%%
In this section we  recall some properties of the fundamental solution $G_\theta=G_\theta(x,t)$ and prove Theorems~\ref{Theorem:1.1}, \ref{Theorem:1.2} and \ref{Theorem:1.3}. 
\vspace{3pt}

Let $0<\theta<2$.
The fundamental solution $G_\theta=G_\theta(x,t)$ is represented by 
\begin{equation}
\label{eq:2.1}
G_\theta(x,t)=(2\pi)^{-\frac{N}{2}}\int_{{\bf R}^N}e^{ix\cdot \xi}e^{-t|\xi|^\theta}\,d\xi
\end{equation}
and it has the following properties: 
\begin{itemize}
  \item[({\bf G})]
  $G_\theta=G_\theta(x,t)$ is a positive smooth function in ${\bf R}^N\times(0,\infty)$ such that 
  \begin{itemize}
  	\item[(i)] 
	$\displaystyle{G_\theta(x,t)=t^{-\frac{N}{\theta}}G_\theta(t^{-\frac{1}{\theta}}x,1)}$ 
	for $x\in{\bf R}^N$ and $t>0$;
	\item[(ii)]
	$\displaystyle{\sup_{x\in{\bf R}^N}\,(1+|x|)^{N+\theta+|\alpha|}|(\partial_x^\alpha G_\theta)(x,1)|<\infty}$ 
	for $\alpha\in{\bf M}$;
	\item[(iii)]
	$G_\theta(\cdot,1)$ is radially symmetric, monotone decreasing with respect to $r:=|x|$ and  
	$$
	\liminf_{|x|\to+\infty}\,(1+|x|)^{N+\theta+j}(\partial_r^j G_\theta)(x,1)>0
	$$
	for $j\in\{0,1,2,\dots\}$;
	\item[(iv)] 
	$\displaystyle{
	G_\theta(x,t)=\int_{{\bf R}^N}G_\theta(x-y,t-s)\,G_\theta(y,s)\,dy}
	$\quad
	for $x\in{\bf R}^N$ and $t>s>0$.
  \end{itemize}
\end{itemize}
See \cite{BT} and \cite{BK}. (See also \cite{IKK01}, \cite{IKK02} and \cite{S}.) 
By ({\bf G})-(i) and (ii), 
for any $\alpha\in{\bf M}$, 
we can find a positive constant $C_\alpha$ such that 
\begin{equation}
\label{eq:2.2}
|(\partial_x^\alpha G_\theta)(x,t)|\le C_\alpha t^{-\frac{N+|\alpha|}{\theta}}
\big(1+t^{-\frac{1}{\theta}}|x|\big)^{-N-\theta-|\alpha|}
\end{equation}
for $x\in{\bf R}^N$ and $t>0$.
This implies that 
\begin{equation}
\label{eq:2.3}	
\sup_{t>0}\,\left[t^{\frac{N}{\theta}(1-\frac{1}{q})+\frac{|\alpha|}{\theta}}\|(\partial_x^\alpha G_\theta)(t)\|_q
+t^{\frac{|\alpha|-\ell}{\theta}}|||(\partial_x^\alpha G_\theta)(t)|||_\ell\right]<\infty
\end{equation}
for $1\le q\le\infty$, $\alpha\in{\bf M}$ and $\ell\in[0,\theta+|\alpha|)$. 
Furthermore, for any $j\in\{0,1,2,\dots\}$, 
by the Young inequality and \eqref{eq:2.3}
we can find a positive constant $C_j$ such that 
\begin{equation}
\label{eq:2.4}
\|\nabla^j[S_\theta(t)\varphi]\|_r\le C_j
t^{-\frac{N}{\theta}(\frac{1}{q}-\frac{1}{r})-\frac{j}{\theta}}\|\varphi\|_q,
\qquad t>0,
\end{equation}
for $\varphi\in L^q$ and $1\le q\le r\le\infty$. 

Next we state a lemma on pointwise estimates of $(\partial_t^m\partial_x^\alpha G_\theta)(x,t)$, 
where $\alpha\in{\bf M}$ and $m\in\{1,2,\dots\}$. 
In what follows,
by the letter $C$
we denote generic positive constants (independent of $x$ and $t$)
and they may have different values also within the same line. 
\begin{lemma}
\label{Lemma:2.1}
Let $\alpha\in{\bf M}$ and $m\in\{1,2,\dots\}$. 
There exists a positive constant $C$ such that 
\begin{equation}
\label{eq:2.5}
|(\partial_t^m\partial_x^\alpha G_\theta)(x,t)|\le Ct^{-\frac{N+|\alpha|}{\theta}-m}
\big(1+t^{-\frac{1}{\theta}}|x|\big)^{-N-\theta m-|\alpha|},
\quad
x\in{\bf R}^N\,\,\,t>0. 
\end{equation}
\end{lemma}
{\bf Proof.}
It follows from \eqref{eq:2.1} that
\begin{equation}
\label{eq:2.6}
|(\partial_t^m\partial_x^\alpha G_\theta)(x,t)|
=t^{-\frac{N+|\alpha|}{\theta}-m}\left|\mathcal{F}^{-1}[f](t^{-\frac{1}{\theta}}x)\right|, 
\qquad 
x\in{\bf R}^N,\,\,t>0,
\end{equation}
where $f(\xi):=\xi^\alpha |\xi|^{\theta m} e^{-|\xi|^\theta}$. 
Then we have
$$
\partial_\xi^{\beta} f\in L^1({\bf R}_\xi^N)\quad\mbox{for all $\beta\in{\bf M}_\Lambda$},
$$
$$
|\partial_\xi^{\beta} f(\xi)|\le C |\xi|^{\theta m+|\alpha|-|\beta|}
\quad \mbox{for all $\xi\in{\bf R}^N\setminus\{0\}$ and $\beta\in{\bf M}_{\Lambda+1}$},
$$
where $\Lambda:=|\alpha|+[\theta m]+N-1\in\{0,1,2,\dots\}$.
By the H\"{o}rmander-Mikhlin type multiplier theorem~(see \cite{SS}) we obtain 
\begin{equation}
\label{eq:2.7}
\left|\mathcal{F}^{-1}[f](x)\right|\le C|x|^{-N-\theta m-|\alpha|},
\qquad x\in{\bf R}^N.
\end{equation} 
On the other hand, it follows that 
\begin{equation}
\label{eq:2.8}
\left|\mathcal{F}^{-1}[f](x)\right|\le C\|f\|_1\le C,
\qquad x\in{\bf R}^N.
\end{equation}
We deduce from \eqref{eq:2.7} and \eqref{eq:2.8} that
$$
\left|\mathcal{F}^{-1}[f](x)\right|\le C(1+|x|)^{-N-\theta m-|\alpha|}, 
\qquad x\in{\bf R}^N.
$$
This together with \eqref{eq:2.6} implies \eqref{eq:2.5}, 
and Lemma~\ref{Lemma:2.1} follows.
\vspace{5pt}
$\Box$

We prepare the following lemma
for the proof of Theorem~\ref{Theorem:1.1}. 
\begin{lemma}
\label{Lemma:2.2}
Let $j\in\{0,1,2,\dots\}$ and $\ell\ge0$. For any $x$, $y\in{\bf R}^N$ and $t>0$, set
$$
H^j_\ell(x,y,t):=(\nabla^j G_\theta)(x-y,t) - \sum_{|\alpha| \le [\ell]}\frac{(-1)^{|\alpha|}}{\alpha !} (\partial^\alpha_x \nabla^j G_\theta)(x,t)y^\alpha. 
$$
\begin{itemize}
  \item[{\rm (a)}]
  There exists $C_1>0$ such that
  \begin{equation}
  \label{eq:2.9}
  \int_{{\bf R}^N} |x|^\ell 
  |H^j_\ell(x,y,t)|\,dx
  \le C_1t^{-\frac{j}{\theta}}|y|^\ell
  \end{equation}
  for $y\in{\bf R}^N$ and $t>0$. 
  \item[{\rm (b)}] 
  There exists $C_2>0$ such that 
  $$
  \int_{{\bf R}^N}
   |x|^\ell\biggr(
   \int_{{\bf R}^N}
   |H^j_K(x,y,t)||\varphi(y)|\,dy\biggr)\,dx
   \le C_2t^{-\frac{K+j-\ell}{\theta}}|||\varphi|||_K,\quad t>0,
  $$
  for $\varphi\in L_K$ and $0\le\ell\le K$. 
  \item[{\rm (c)}] 
  For any $\varphi\in L_K$ and $0\le\ell\le K$, 
  \begin{equation*}
  \lim_{t\to\infty}t^{\frac{K+j-\ell}{\theta}}
  \int_{{\bf R}^N}
   |x|^\ell\biggr|
   \int_{{\bf R}^N}
   H^j_K(x,y,t)\varphi(y)\,dy\biggr|\,dx=0. 
  \end{equation*}
  \end{itemize}

  \end{lemma}
{\bf Proof.} 
We prove assertion~(a). Let $y\in{\bf R}^N$. 
Since 
\begin{equation}
\label{eq:2.10}
\begin{split}
H^j_\ell(x,y,t)
 & =\frac{1}{[\ell]!}\int_0^1(1-\tau)^{[\ell]}\frac{d^{[\ell]+1}}{d\tau^{[\ell]+1}}(\nabla^j G_\theta)(x-\tau y,t)\,d\tau\\
 & =(-1)^{[\ell]+1}([\ell]+1)\sum_{|\alpha|=[\ell]+1}\frac{y^\alpha}{\alpha!}
 \int_0^1(1-\tau)^{[\ell]}(\partial_x^\alpha \nabla^j G_\theta)(x-\tau y,t)\,d\tau,
\end{split}
\end{equation}
by \eqref{eq:2.2}
%({\bf G})-(ii) 
we have
\begin{equation}
\label{eq:2.11}
\begin{split}
& \int_{\{|x| \ge 2|y|\}} |x|^\ell 
|H^j_\ell(x,y,t)|\,dx \\
& \le C\int_0^1\int_{\{|x| \ge 2|y|\}} |x|^\ell \left|(\nabla^{[\ell]+j+1} G_\theta)(x-\tau y,t) \right| |y|^{[\ell]+1}\,dx\,d\tau \\
& \le C|y|^\ell \int_0^1\int_{\{|x| \ge 2|y|\}} |x|^{[\ell]+1} \left|(\nabla^{[\ell]+j+1} G_\theta)(x-\tau y,t) \right|\,dx\,d\tau \\
&\le C|y|^\ell \int_0^1\int_{\{|x| \ge 2|y|\}} |x|^{[\ell]+1} \ t^{-\frac{N}{\theta}-\frac{[\ell]+j+1}{\theta}} 
\left(1+t^{-\frac{1}{\theta}}|x-\tau y| \right)^{-(N+\theta+[\ell]+j+1)}\,dx\,d\tau. 
\end{split}
\end{equation}
It follows that 
$$
|x-\tau y| \ge |x|-|y| \ge |x|/2
\quad\mbox{if}\quad
|x| \ge 2|y|\quad\mbox{and}\quad
0 \le \tau \le 1.
$$
This together with \eqref{eq:2.11} implies that 
\begin{equation}
\label{eq:2.12}
\begin{split}
& \int_{\{|x| \ge 2|y|\}} |x|^\ell 
|H^j_\ell(x,y,t)|\,dx \\
&\le C|y|^\ell \int_{{\bf R}^N} |x|^{[\ell]+1} \ t^{-\frac{N}{\theta}-\frac{[\ell]+j+1}{\theta}} 
\left(1+t^{-\frac{1}{\theta}}\frac{|x|}{2} \right)^{-(N+\theta+[\ell]+j+1)}\,dx\le Ct^{-\frac{j}{\theta}}|y|^\ell.
\end{split}
\end{equation}
On the other hand, by \eqref{eq:2.2} we see that 
\begin{equation}
\label{eq:2.13}
\begin{split}
& \int_{\{|x| < 2|y|\}} |x|^\ell 
|H^j_\ell(x,y,t)|\,dx\\
 & \le \int_{\{|x| < 2|y|\}} |x|^\ell |(\nabla^j G_\theta)(x-y,t)|\,dx\\
 & \qquad\quad
+C\sum_{|\alpha|\le[\ell]} \int_{\{|x| < 2|y|\}} |x|^\ell |(\partial_x^\alpha\nabla^j G_\theta)(x,t)||y|^{|\alpha|}\,dx \\
 & \le (2|y|)^\ell\int_{{\bf R}^N}|(\nabla^j G_\theta)(x,t)|\,dx+C\sum_{|\alpha|\le[\ell]}  |y|^\ell  \int_{\{|x| < 2|y|\}} |x|^{|\alpha|} |(\partial_x^\alpha\nabla^j G_\theta)(x,t)|\,dx\\
 & \le Ct^{-\frac{j}{\theta}}|y|^\ell.
\end{split}
\end{equation}
By \eqref{eq:2.12} and \eqref{eq:2.13} we obtain \eqref{eq:2.9}. 
Thus assertion~(a) follows. 

We prove assertions~(b) and (c). 
It follows from Lemma~\ref{Lemma:2.2}~(a) and \eqref{eq:2.2} that
\begin{equation}
\label{eq:2.14}
\begin{split}
 & \int_{\{|y|\ge R^{\frac{1}{\theta}}\}}
 \biggr(
 \int_{{\bf R}^N}
 |x|^\ell |H^j_K(x,y,t)|\,dx
 \biggr)|\varphi(y)|\,dy\\
  & \le\int_{\{|y|\ge R^{\frac{1}{\theta}}\}}
 \biggr(
 \int_{{\bf R}^N}
 |x|^\ell |H^j_\ell(x,y,t)|\,dx
 \biggr)|\varphi(y)|\,dy\\
 & \qquad\qquad
 +C\sum_{[\ell]<|\alpha|\le K}
 \int_{\{|y|\ge R^{\frac{1}{\theta}}\}}
 \biggr(\int_{{\bf R}^N}|x|^\ell|(\partial_x^\alpha \nabla^j G_\theta)(x,t)|\,dx\biggr)|y|^{|\alpha|}|\varphi(y)|\,dy\\ 
 & \le Ct^{-\frac{j}{\theta}}\int_{\{|y|\ge R^{\frac{1}{\theta}}\}}|y|^\ell|\varphi(y)|\,dy
 +C\sum_{[\ell]<|\alpha|\le K}t^{-\frac{|\alpha|+j-\ell}{\theta}}\int_{\{|y|\ge R^{\frac{1}{\theta}}\}}|y|^{|\alpha|}|\varphi(y)|\,dy\\
 & \le Ct^{-\frac{j}{\theta}}\int_{\{|y|\ge R^{\frac{1}{\theta}}\}}|y|^\ell\left(\frac{|y|}{R^{\frac{1}{\theta}}}\right)^{K-\ell}|\varphi(y)|\,dy\\
 & \qquad\qquad
 +C\sum_{[\ell]<|\alpha|\le K}t^{-\frac{|\alpha|+j-\ell}{\theta}}
 \int_{\{|y|\ge R^{\frac{1}{\theta}}\}}\left(\frac{|y|}{R^{\frac{1}{\theta}}}\right)^{K-|\alpha|}|y|^{|\alpha|}|\varphi(y)|\,dy\\
 & =Ct^{-\frac{K+j-\ell}{\theta}}
 \biggr[(R^{-1}t)^{\frac{K-\ell}{\theta}}+\sum_{[\ell]<|\alpha|\le K}(R^{-1}t)^{\frac{K-|\alpha|}{\theta}}\biggr]
 \int_{\{|y|\ge R^{\frac{1}{\theta}}\}}|y|^K|\varphi(y)|\,dy
\end{split}
\end{equation}
for all $R>0$. 
On the other hand, by \eqref{eq:2.2} and \eqref{eq:2.10} we obtain 
\begin{equation}
\label{eq:2.15}
\begin{split}
 & \int_{\{|y|<R^{\frac{1}{\theta}}\}}
 \biggr(
 \int_{{\bf R}^N}
 |x|^\ell |H^j_K(x,y,t)|\,dx
 \biggr)|\varphi(y)|\,dy\\
 & \le C\int_0^1\int_{\{|y|<R^{\frac{1}{\theta}}\}} 
 \biggr(\int_{{\bf R}^N} |x|^\ell |(\nabla^{[K]+j+1}G_\theta)(x-\tau y,t)| |y|^{[K]+1}\,dx\biggr)|\varphi(y)|\,dy\,d\tau \\
 & = C\int_0^1\int_{\{|y|<R^{\frac{1}{\theta}}\}} 
 \biggr(\int_{{\bf R}^N} |x+\tau y|^\ell |(\nabla^{[K]+j+1}G_\theta)(x,t)|\,dx\biggr)|y|^{[K]+1}|\varphi(y)|\,dy \,d\tau \\
 & \le C\int_{\{|y|<R^{\frac{1}{\theta}}\}}
 \biggr(\int_{{\bf R}^N} (|x|^\ell +|y|^\ell) |(\nabla^{[K]+j+1}G_\theta)(x,t)|\,dx\biggr)|y|^{[K]+1} |\varphi(y)|\,dy\\
 & \le C\int_{\{|y|<R^{\frac{1}{\theta}}\}} 
 (t^{-\frac{[K]+j+1-\ell}{\theta}}+t^{-\frac{[K]+j+1}{\theta}} |y|^{\ell})|y|^{[K]+1} |\varphi(y)|\,dy\\
 & \le Ct^{-\frac{j}{\theta}}(t^{-\frac{[K]+1-\ell}{\theta}}R^{\frac{[K]+1-K}{\theta}}
 +t^{-\frac{[K]+1}{\theta}}R^{\frac{[K]+\ell+1-K}{\theta}})
 \int_{\{|y|<R^{\frac{1}{\theta}}\}}|y|^K |\varphi(y)|\,dy
\end{split}
\end{equation}
for all $R>0$. 
Then, by \eqref{eq:2.14} and \eqref{eq:2.15}
we set $R=t$ to obtain 
$$
\int_{{\bf R}^N}|x|^\ell\biggr(
\int_{{\bf R}^N}|H^j_K(x,y,t)||\varphi(y)|\,dy\biggr)\,dx
\le Ct^{-\frac{K+j-\ell}{\theta}}|||\varphi|||_K,\qquad t>0.
$$
This implies assertion~(b). 
Similarly, setting $R=\epsilon t$ with $0<\epsilon\le 1$, 
we have 
\begin{equation*}
\begin{split}
 & \int_{{\bf R}^N}|x|^\ell\biggr|
   \int_{{\bf R}^N}H^j_K(x,y,t)
   \varphi(y)\,dy\biggr|\,dx\\
 & \le Ct^{-\frac{K+j-\ell}{\theta}}
\biggr[(\epsilon^{-1})^{\frac{K-\ell}{\theta}}+\sum_{[\ell]<|\alpha|\le K}(\epsilon^{-1})^{\frac{K-|\alpha|}{\theta}}\biggr]
\int_{\{|y| \ge(\epsilon t)^{\frac{1}{\theta}}\}}|y|^K|\varphi(y)|\,dy\\
 & +Ct^{-\frac{K+j-\ell}{\theta}}(\epsilon^{\frac{[K]+1-K}{\theta}}+\epsilon^{\frac{[K]+\ell+1-K}{\theta}})
|||\varphi|||_K.
\end{split}
\end{equation*}
This implies that 
\begin{equation*}
\begin{split}
 & \limsup_{t\to\infty}\,t^{\frac{K+j-\ell}{\theta}}
 \int_{{\bf R}^N}|x|^\ell\biggr|
  \int_{{\bf R}^N}H^j_K(x,y,t)
  \varphi(y)\,dy\biggr|\,dx\\
 & \le C(\epsilon^{\frac{[K]+1-K}{\theta}}+\epsilon^{\frac{[K]+\ell+1-K}{\theta}})
|||\varphi|||_K. 
\end{split}
\end{equation*}
Then, since $\epsilon$ is arbitrary,  
we obtain assertion~(c). 
Thus Lemma~\ref{Lemma:2.2} follows. 
$\Box$
\vspace{5pt}

Now we are ready to prove Theorem~\ref{Theorem:1.1}.
\vspace{3pt}
\newline
{\bf Proof of Theorem~\ref{Theorem:1.1}.}
Let $\varphi\in L_K$ and $j\in\{0,1,2,\dots\}$. 
It follows from \eqref{eq:1.8} and \eqref{eq:2.10} that 
\begin{equation}
\label{eq:2.16}
\begin{split}
(\nabla^j v)(x,t)
 & =\nabla^j [S_\theta(t)\varphi](x)
-\sum_{|\alpha| \le K}\biggr[\int_{{\bf R}^N}y^\alpha \varphi(y)\,dy\biggr](\nabla^j g_{\alpha,\theta})(x,t)\\
 & =\int_{{\bf R}^N}
\biggr[(\nabla^j G_\theta)(x-y,t)-\sum_{|\alpha| \le K}\frac{(-1)^{|\alpha|}}{\alpha !} 
(\partial^\alpha_x \nabla^j G_\theta)(x,t)y^\alpha\biggr]\varphi(y)\,dy\\
 & =\int_{{\bf R}^N} H^j_K(x,y,t)\varphi(y)\,dy
\end{split}
\end{equation}
for $x\in{\bf R}^N$ and $t>0$. 
Since
$$
(\nabla^j v)(x,t)=\biggr[S_\theta\biggr(\frac{t}{2}\biggr)(\nabla^j v)\biggr(\cdot,\frac{t}{2}\biggr)\biggr](x),
$$ 
by \eqref{eq:2.4} we have
$$
t^{\frac{N}{\theta}(1-\frac{1}{q})}\|(\nabla^j v)(t)\|_q+t^{-\frac{\ell}{\theta}}|||(\nabla^j v)(t)|||_\ell 
\le C\|(\nabla^j v)(t/2)\|_1+t^{-\frac{\ell}{\theta}}|||(\nabla^j v)(t)|||_\ell, \quad t>0.
$$
Then Lemma~\ref{Lemma:2.2}~(b) and (c) with \eqref{eq:2.16} imply assertions~(a) and (b), respectively. 

It remans to prove \eqref{eq:1.9}. 
Let $0\le|\alpha|\le K$. 
By Lemma~\ref{Lemma:2.2}~(b) and \eqref{eq:2.16} 
we apply the Fubini theorem to obtain 
\begin{equation}
\label{eq:2.17}
\begin{split}
\int_{{\bf R}^N} x^\alpha (\nabla^j v)(x,t)\,dx
 & =\int_{{\bf R}^N}
 x^\alpha\biggr[\int_{{\bf R}^N}
H^j_K(x,y,t)\varphi(y)\,dy\biggr]\,dx\\
& =\int_{{\bf R}^N}\biggr(\int_{{\bf R}^N}x^\alpha H^j_K(x,y,t)\,dx\biggr)\varphi(y)\,dy
\end{split}
\end{equation}
for $t>0$. 
On the other hand, it follows that 
\begin{equation}
\label{eq:2.18}
\int_{{\bf R}^N}x^\alpha(\partial_x^\beta G_\theta)(x,t)\,dx=0 \quad \mbox{if not} \quad \alpha\ge\beta,
\end{equation}
in the case $|\alpha|<|\beta|+\theta$. 
Then, by \eqref{eq:2.10} and \eqref{eq:2.18} we have 
\begin{equation*}
\begin{split}
 & \int_{{\bf R}^N}x^\alpha H^j_K(x,y,t)\,dx\\
 & =(-1)^{[K]+1}([K]+1)\sum_{|\beta|=[K]+1}\frac{1}{\beta!}\int_{{\bf R}^N}x^\alpha
\int_0^1 (1-\tau)^{[K]}(\partial_x^\beta \nabla^j G_\theta)(x-\tau y,t)y^{\beta}\,d\tau\,dx\\
 & =(-1)^{[K]+1}([K]+1)\sum_{|\beta|=[K]+1}\frac{y^\beta}{\beta!}
 \int_0^1(1-\tau)^{[K]}\biggr(\int_{{\bf R}^N}x^\alpha(\partial_x^\beta \nabla^j G_\theta)(x-\tau y,t)\,dx\biggr)\,d\tau=0
\end{split}
\end{equation*}
for $y\in{\bf R}^N$ and $t>0$. 
This together with \eqref{eq:2.17} implies that 
$$
\int_{{\bf R}^N} x^\alpha(\nabla^jv)(x,t)\,dx=0,\qquad t>0.
$$
Thus \eqref{eq:1.9} holds and the proof of Theorem~\ref{Theorem:1.1} is complete. 
$\Box$
\vspace{8pt}
\newline
By Theorem~\ref{Theorem:1.1} we prove Theorems~\ref{Theorem:1.2} and \ref{Theorem:1.3}.
\vspace{3pt}
\newline
{\bf Proof of Theorem~\ref{Theorem:1.2}.}
Let $j\in\{0,1,2,\dots\}$ and $\varphi\in \mathcal{L}_{K,j}^0$. 
By the definition of $\mathcal{L}_{K,j}^0$ if $0\le K<\theta+j$, then we see that $\mathcal{L}_{K,j}^0$ and the assertion~(a) holds. So it suffices to consider the case $K\ge\theta+j$. By Theorem~\ref{Theorem:1.1} we have 
$$
\nabla^j[S_\theta(t)\varphi]
-\sum_{K-j-\theta<|\alpha| \le K}\biggr[\int_{{\bf R}^N}
y^\alpha \varphi(y)\,dy\biggr](\nabla^j g_{\alpha,\theta})(t)\in{\mathcal L}_{K,j}^0,
\qquad t>0.
$$
If $|\alpha|>K-j-\theta$, then it follows from \eqref{eq:2.3} and \eqref{eq:2.18} 
that $(\nabla^j g_{\alpha,\theta})(t)\in{\mathcal L}_{K,j}^0$ for $t>0$. 
Therefore ${\mathcal L}_{K,j}^0$ satisfies \eqref{eq:1.11}. 
On the other hand, 
if $\varphi\in L_K\setminus{\mathcal L}_{K,j}^0$, then it follows from Theorem~\ref{Theorem:1.1} and ({\bf G})-(ii), (iii) 
that $\nabla^j [S_\theta(t)\varphi]\not\in L_K$ for any $t>0$. 
Thus assertion~(a) follows. 

We prove assertion~(b). 
By Theorem~\ref{Theorem:1.1} we see that 
\begin{equation}
\label{eq:2.19}
\begin{split}
\left\{\nabla^j [S_\theta(t)\varphi]\,:\,\varphi\in{\mathcal L}_{K,j}\right\}
\subset L_K
& +\sum_{|\alpha|\le K}\langle(\partial_x^\alpha \nabla^j G_\theta)(\cdot,t)\rangle \\
& +\sum_{|\alpha|+m\theta\le K-j}\langle(\partial_t^m\partial_x^\alpha \nabla^j G_\theta)(\cdot,1+t)\rangle
\end{split}
\end{equation}
for $t>0$. 
For any $\alpha\in{\bf M}$ and $m\in\{0,1,2,\dots\}$ 
with $|\alpha|+\theta m\le K-j$, 
let $m'$ be the largest nonnegative integer satisfying $|\alpha|+\theta(m+m')\le K-j$. 
Since
\begin{equation*}
\begin{split}
 &  (\partial_t^m\partial_x^\alpha \nabla^j G_\theta)(x,t)
=\sum_{i=0}^{m'}\frac{(t-1)^i}{i!}(\partial_t^{m+i}\partial_x^\alpha \nabla^j G_\theta)(x,1) \\
 & \qquad\quad
+\frac{(t-1)^{m'+1}}{(m')!}\int_0^1(1-\tau)^{m'}
\frac{d^{m'+1}}{d\tau^{m'+1}}[(\partial_t^m\partial_x^\alpha\nabla^j  G_\theta)(x,1+\tau(t-1))]\,d\tau,
\end{split}
\end{equation*}
by Lemma~\ref{Lemma:2.1} and ({\bf G})-(ii) we have
$$
(\partial_t^m\partial_x^\alpha\nabla^j G_\theta)(\cdot,t)
\in\sum_{i=0}^{m'}\langle (\partial_t^{m+i}\partial_x^\alpha \nabla^j G_\theta)(\cdot,1)\rangle
+L_K\subset\mathcal{L}_{K,j}\subset L_{K'},\qquad t>0,
$$
for any $0\le K'\le K$ with $K'<\theta$. 
This together with \eqref{eq:2.19} implies assertion~(b). 
Thus Theorem~\ref{Theorem:1.2} follows. 
$\Box$
\vspace{3pt}
\newline
{\bf Proof of Theorem~\ref{Theorem:1.3}.}
We can assume, without loss of generality, that $C(\varphi)=0$, 
that is
\begin{equation}
\label{eq:2.20}
\int_{{\bf R}^N}x\varphi(x)\,dx=0. 
\end{equation}
We prove that $H(t)\not=\emptyset$ for any $t>0$. 
For any $t>0$, since 
$$
\int_{{\bf R}^N}[S_\theta(t)\varphi](x)\,dx=\int_{{\bf R}^N}\varphi(x)\,dx>0, 
$$
we can find $x_t\in{\bf R}^N$ such that $[S_\theta(t)\varphi](x_t)>0$.  
By \eqref{eq:2.2} 
we can find $L>0$ such that 
\begin{equation*}
\begin{split}
[S_\theta(t)\varphi](x)
 & \le Ct^{-\frac{N}{\theta}}\int_{{\bf R}^N}(1+t^{-\frac{1}{\theta}}|x-y|)^{-N-\theta}|\varphi(y)|\,dy\\
 & \le Ct^{-\frac{N}{\theta}}\int_{B(0,L)}(1+t^{-\frac{1}{\theta}}|x-y|)^{-N-\theta}|\varphi(y)|\,dy
+Ct^{-\frac{N}{\theta}}\int_{{\bf R}^N\setminus B(0,L)}|\varphi(y)|\,dy\\
 & \le Ct^{-\frac{N}{\theta}}(1+t^{-\frac{1}{\theta}}L)^{-N-\theta}\|\varphi\|_1
 +Ct^{-\frac{N}{\theta}}\int_{{\bf R}^N\setminus B(0,L)}|\varphi(y)|\,dy\\
 & <[S_\theta(t)\varphi](x_t)
\end{split}
\end{equation*}
for $x\in{\bf R}^N\setminus B(0,2L)$. 
This means that $H(t)\not=\emptyset$ for any $t>0$. 

We study the behavior of $H(t)$ as $t\to\infty$. 
It follows from \eqref{eq:2.20} and Theorem~\ref{Theorem:1.1} with $q=\infty$, $j=0$ and $K=1$ that 
\begin{equation}
\label{eq:2.21}
\left\|S_\theta(t)\varphi-M(\varphi)g_{0,\theta}(t)\right\|_{\infty}
\le Ct^{-\frac{N+1}{\theta}},\qquad t>0.
\end{equation} 
For any $\epsilon>0$, 
by ({\bf G})-(iii) and \eqref{eq:2.21} we have  
\begin{equation*}
\begin{split}
[S_\theta(t)\varphi](x)-[S_\theta(t)\varphi](0)
 & =M(\varphi)t^{-\frac{N}{\theta}}\left\{G_\theta(t^{-\frac{1}{\theta}}x,1)-G_\theta(0,1)\right\}+O(t^{-\frac{N+1}{\theta}})\\
 & \le -M(\varphi)t^{-\frac{N}{\theta}}\{G_\theta(0,1)-G_\theta(\epsilon,1)\}+O(t^{-\frac{N+1}{\theta}})<0
\end{split}
\end{equation*}
for all sufficiently large $t>0$ 
uniformly on the set $\{x\in{\bf R}^N\,:|x|\ge\epsilon t^{\frac{1}{\theta}}\}$.
This implies that 
\begin{equation}
\label{eq:2.22}
H(t)\subset B(0,\epsilon t^{\frac{1}{\theta}})
\end{equation}
for all sufficiently large $t>0$. 

By \eqref{eq:2.20} we apply Theorem~\ref{Theorem:1.1} with $q=\infty$, $j=0$ and $K=2$ to have
\begin{equation}
\label{eq:2.23}
[S_\theta(t)\varphi](x)
=M(\varphi)G_\theta(x,t)+\sum_{|\alpha|=2}c_\alpha\,g_{\alpha,\theta}(x,t)
+o\left(t^{-\frac{N+2}{\theta}}\right)
\end{equation}
as $t\to\infty$ uniformly for $x\in{\bf R}^N$, where 
$$
c_\alpha:=\int_{{\bf R}^N}x^\alpha\varphi(x)\,dx.
$$
Taking a sufficiently small $\epsilon>0$ if necessary, 
by \eqref{eq:1.7}, \eqref{eq:2.20} and \eqref{eq:2.23} 
we obtain 
\begin{equation}
\label{eq:2.24}
\begin{split}
 & [S_\theta(t)\varphi](x)\\
 & =M(\varphi)t^{-\frac{N}{\theta}}G_\theta(t^{-\frac{1}{\theta}}x,1)
+t^{-\frac{N+2}{\theta}}\sum_{|\alpha|=2}\frac{c_\alpha}{\alpha!}
(\partial_x^\alpha G_\theta)(t^{-\frac{1}{\theta}}x,1)+o\left(t^{-\frac{N+2}{\theta}}\right)\\
 & =M(\varphi)t^{-\frac{N}{\theta}}
 \biggr[G_\theta(0,1)+t^{-\frac{2}{\theta}}\sum_{|\alpha|=2}\frac{1}{\alpha!}(\partial_x^\alpha G_\theta)(0,1)x^\alpha
 +O\left(t^{-\frac{3}{\theta}}|x|^3\right)\biggr]\\
 & \qquad\qquad\qquad
 +t^{-\frac{N+2}{\theta}}\sum_{|\alpha|=2}\frac{c_\alpha}{\alpha!}
\left[(\partial_x^\alpha G_\theta)(0,1)+O(t^{-\frac{1}{\theta}}|x|)\right]+o\left(t^{-\frac{N+2}{\theta}}\right)\\
\end{split}
\end{equation}
for all sufficiently large $t>0$ uniformly on $\{x\in{\bf R}^N\,:\,|x|\le\epsilon t^{\frac{1}{\theta}}\}$. 
On the other hand, by \eqref{eq:2.1} we have 
\begin{equation}
\label{eq:2.25}
\frac{\partial^2 G_\theta}{\partial x_i^2}(0,1)=-(2\pi)^{-\frac{N}{2}}\int_{{\bf R}^N}\xi_i^2 e^{-|\xi|^\theta}\,d\xi<0,
\qquad
i=1,\dots,N.
\end{equation}
Taking a sufficiently small $\epsilon>0$ if necessary again, by \eqref{eq:2.24} and \eqref{eq:2.25}, 
for any $R>0$, we obtain 
\begin{equation}
\label{eq:2.26}
\begin{split}
 & [S_\theta(t)\varphi](x)-[S_\theta(t)\varphi](0)\\
 & =M(\varphi)t^{-\frac{N+2}{\theta}}
 \sum_{|\alpha|=2}\frac{1}{\alpha!}
 (\partial_x^\alpha G_\theta)(0,1)x^\alpha+O(t^{-\frac{N+3}{\theta}}(|x|+|x|^3))
 +o\left(t^{-\frac{N+2}{\theta}}\right)\\
 & \le -CM(\varphi)t^{-\frac{N+2}{\theta}}|x|^2
 +O\left(\epsilon t^{-\frac{N+2}{\theta}}\right)
 +O\left(\epsilon t^{-\frac{N+2}{\theta}}|x|^2\right)
 +o\left(t^{-\frac{N+2}{\theta}}\right)\\
  & \le -CM(\varphi)t^{-\frac{N+2}{\theta}}|x|^2
 +O\left(\epsilon t^{-\frac{N+2}{\theta}}\right)
 +o\left(t^{-\frac{N+2}{\theta}}\right)\\
 & \le -CM(\varphi)t^{-\frac{N+2}{\theta}}R^2
 +O\left(\epsilon t^{-\frac{N+2}{\theta}}\right)
 +o\left(t^{-\frac{N+2}{\theta}}\right)
\end{split}
\end{equation}
for $x\in{\bf R}^N$ with $R\le|x|\le\epsilon t^{\frac{1}{\theta}}$ and sufficiently large $t>0$. 
Since $R$ is arbitrary, 
we deduce from \eqref{eq:2.22} and \eqref{eq:2.26} that 
\begin{equation}
\label{eq:2.27}
\lim_{t\to\infty}\sup_{x\in H(t)}\,|x|=0.
\end{equation}
Furthermore, by \eqref{eq:2.1} and \eqref{eq:2.25} 
we apply Theorem~\ref{Theorem:1.1} with $q=\infty$, $j=2$ and $K=2$ to obtain
\begin{equation}
\label{eq:2.28}
\begin{split}
\left(\nabla^2[S_\theta(t)\varphi](x)y,y\right)
 & =M(\varphi)t^{-\frac{N+2}{\theta}}((\nabla^2 G_\theta)(0,1)y,y)+o(t^{-\frac{N+2}{\theta}})|y|^2\\
 & \le -CM(\varphi)t^{-\frac{N+2}{\theta}}|y|^2,
\quad y\in{\bf R}^N,
\end{split}
\end{equation}
for $x\in{\bf R}^N$ with $|x|\le\epsilon t^{\frac{1}{\theta}}$ and sufficiently large $t>0$. 
By \eqref{eq:2.27} and \eqref{eq:2.28} we see that 
$H(t)$ consists of one point $x(t)\in{\bf R}^N$ for all sufficiently large $t>0$. 
In addition, the implicit function theorem implies that 
$\{x(t)\}$ moves along a smooth curve in ${\bf R}^N$. 
Thus Theorem~\ref{Theorem:1.3} follows. 
$\Box$
%%%%%%%%%%%%%%%%%%%%%%%%%%%%%%%%%%%%%
%%%%%%%%%%%%%%%%%%%%%%%%%%%%%%%%%%%%%
\section{Fractional diffusion equation with inhomogeneous term}
%%%%%%%%%%%%%%%%%%%%%%%%%%%%%%%%%%%%%
%%%%%%%%%%%%%%%%%%%%%%%%%%%%%%%%%%%%%
Let $K\ge 0$. 
Let $\psi$ be a radially symmetric smooth function in ${\bf R}^N$ such that 
\begin{equation}
\label{eq:3.1}
\partial_x^\alpha\psi\in L^\infty\cap L_{K+|\alpha|}\,\,\,\mbox{for $\alpha\in{\bf M}_K$},
\qquad
\int_{{\bf R}^N}\psi(y)\,dy=1. 
\end{equation} 
For any $\alpha\in{\bf M}$, $x\in{\bf R}^N$, $t>0$ and $s\ge0$, we define 
\begin{equation}
\label{eq:3.2}
\begin{split}
 & \psi(x,s):=(1+s)^{-\frac{N}{\theta}}\psi\left((1+s)^{-\frac{1}{\theta}}x\right),\\
 & \psi_\alpha(x,s):=(1+s)^{-\frac{N+|\alpha|}{\theta}}\frac{(-1)^{|\alpha|}}{\alpha!}
(\partial_x^\alpha\psi)\left((1+s)^{-\frac{1}{\theta}}x\right),\\
 & \Psi_{\alpha,\theta}(x,t:s):=\frac{(-1)^{|\alpha|}}{\alpha!}\partial_x^\alpha [S_\theta(t)\psi(\cdot,s)](x).
\end{split}
\end{equation}
It follows from \eqref{eq:3.1} and \eqref{eq:3.2} that
\begin{equation}
\label{eq:3.3}
(1+s)^{\frac{N}{\theta}(1-\frac{1}{q})+\frac{|\alpha|}{\theta}}\|\psi_\alpha(s)\|_q
+(1+s)^{-\frac{\ell}{\theta}+\frac{|\alpha|}{\theta}}|||\psi_\alpha(s)|||_\ell<+\infty, \qquad s\ge0,
\end{equation}
for $1\le q\le\infty$ and $0\le\ell<K+|\alpha|$. 
Furthermore, 
\begin{equation}
\label{eq:3.4}
\Psi_{\alpha,\theta}(x,t:s)=[S_\theta(t)\psi_\alpha(\cdot,s)](x)
\end{equation}
for $\alpha\in{\bf M}_K$, $x\in{\bf R}^N$, $t>0$ and $s\ge0$. 
The precise description of the asymptotic behavior of $\Psi_{\alpha,\theta}(\cdot,t:s)$ as $t\to\infty$ 
can be given by Theorem~\ref{Theorem:1.1} (see also Theorem~\ref{Theorem:1.2}).  

In this section, 
by using $\{\Psi_{\alpha,\theta}\}$ 
we obtain the asymptotic expansions of 
the solution~$u$ to the fractional diffusion equation with an inhomogeneous term $F=F(x,t)$
\begin{equation}
\label{eq:3.5}
\partial_t u+(-\Delta)^{\frac{\theta}{2}}u=F\quad\mbox{in}\quad{\bf R}^N\times(0,\infty),
\qquad
u(x,0)=\varphi(x)\quad\mbox{in}\quad{\bf R}^N,
\end{equation}
where $\varphi\in L^\infty\cap L_K$ and $F\in L^\infty(0,T:L_K)$ for any $T>0$. 
Similarly to \eqref{eq:1.2}, we say that $u$ is a solution of \eqref{eq:3.5} if 
$u$ is a continuous function in ${\bf R}^N\times(0,\infty)$ 
and $u$ satisfies
\begin{equation*}
\begin{split}
u(x,t) & =\int_{{\bf R}^N}G_\theta(x-y,t)\varphi(y)\,dy+\int_0^t\int_{{\bf R}^N}G_\theta(x-y,t-s)F(y,s)\,dy\,ds\\
 & \equiv [S_\theta(t)\varphi](x)+\int_0^t [S_\theta(t-s)F(\cdot,s)](x)\,ds,
 \qquad x\in{\bf R}^N,\,\,t>0. 
\end{split}
\end{equation*}

We state the main result of this section, which is 
a generalization of \cite[Theorems~1.1 and 1.2]{IKK02} (see Remark~\ref{Remark:3.1}~(i)).
In what follows, 
for any $f\in L_K$, $\alpha\in{\bf M}_K$ and $s\ge 0$,
we define 
$M_\alpha(f,s)$ inductively in $\alpha$ by 
\begin{equation}
\label{eq:3.6}
\left\{
\begin{array}{l}
\displaystyle{M_0(f,s):=\int_{{\bf R}^N}f(x)\,dx}\quad\mbox{if}\quad \alpha=0,\vspace{7pt}\\
\displaystyle{M_\alpha(f,s):=\int_{{\bf R}^N} 
x^\alpha f(x)\,dx-\sum_{\beta\le\alpha,\beta\not=\alpha}M_\beta(f,s) 
\int_{{\bf R}^N}x^\alpha\psi_\beta(x,s)\,dx}\quad\mbox{if}\quad\alpha\not=0.
\end{array}
\right.
\end{equation}
Then it follows from \eqref{eq:3.1} and \eqref{eq:3.6} that 
\begin{equation}
\label{eq:3.7}
\int_{{\bf R}^N}
x^\beta\biggr[f(x)-\sum_{|\alpha|\le K}M_\alpha(f,s)\psi_\alpha(x,s)\biggr]\,dx=0
\end{equation}
for $\beta\in{\bf M}_K$ and $s\ge 0$. 
\begin{theorem}
\label{Theorem:3.1}
Let $N\ge1$, $0<\theta<2$ and $K\ge 0$. 
\begin{itemize}
  \item[{\rm (i)}] 
  Let $\varphi\in L^\infty\cap L_K$ and set 
  \begin{equation*}
  w_1(x,t):=\sum_{|\alpha|\le K}M_\alpha(\varphi,0)\Psi_{\alpha,\theta}(x,t:0).
  \end{equation*} 
  Then there exists a constant $C_1>0$ such that 
  \begin{equation}
  \label{eq:3.8}
  \begin{split}
   & (1+t)^{\frac{N}{\theta}(1-\frac{1}{q})}\|S_\theta(t)\varphi-w_1(t)\|_q
  +(1+t)^{-\frac{\ell}{\theta}}|||S_\theta(t)\varphi-w_1(t)|||_\ell\\
   & \qquad\quad
  \le C_1(1+t)^{-\frac{K}{\theta}}(\|\varphi\|_\infty+\|\varphi\|_1+|||\varphi|||_K)
  \end{split}
  \end{equation}
  for $t>0$, $1\le q\le\infty$ and $0\le \ell\le K$. 
  In particular, for any $1\le q\le\infty$ and $0\le\ell\le K$,
  \begin{equation}
  \label{eq:3.9}
  \lim_{t\to\infty}
  t^{\frac{K}{\theta}}
  \biggr[t^{\frac{N}{\theta}(1-\frac{1}{q})}\|S_\theta(t)\varphi-w_1(t)\|_q
  +t^{-\frac{\ell}{\theta}}|||S_\theta(t)\varphi-w_1(t)|||_\ell\biggr]=0.
  \end{equation}
  \item[{\rm (ii)}] 
  Let $1\le q\le\infty$ and let $F$ be a measurable function in ${\bf R}^N\times(0,\infty)$ such that 
  \begin{equation}
  \label{eq:3.10}
  \begin{split}
  E_{K,q}[F](t):= & (1+t)^{\frac{K}{\theta}}
  \left[(1+t)^{\frac{N}{\theta}(1-\frac{1}{q})}\|F(t)\|_q+\|F(t)\|_1\right]+|||F(t)|||_K\\
  \in & L^\infty(0,T)
  \end{split}
  \end{equation}
  for any $T>0$.  
  Set 
  \begin{equation*}
  \begin{split}
  v(x,t):= & \int_0^t [S_\theta(t-s)F(\cdot,s)](x)\,ds,\\
  w_2(x,t):= & \sum_{|\alpha|\le K}
  \int_0^t M_\alpha(F(s),s)\Psi_{\alpha,\theta}(x,t-s:s)\,ds.
  \end{split}
  \end{equation*}
  Then there exists a constant $C_2>0$ such that
  \begin{equation}
  \label{eq:3.11}
  \|v(t)-w_2(t)\|_q+|||v(t)-w_2(t)|||_\ell
  \le C_2\int_0^t E_{K,q}[F](s)\,ds,\quad t>0,
  \end{equation}
  for $0\le\ell\le K$. 
  Furthermore, there exists a constant $C_3>0$ such that, 
  for any $\epsilon>0$ and $T>0$, 
  \begin{equation}
  \label{eq:3.12}
  \begin{split}
   & t^{\frac{N}{\theta}(1-\frac{1}{q})}\|v(t)-w_2(t)\|_q
  +t^{-\frac{\ell}{\theta}}|||v(t)-w_2(t)|||_\ell\\
   & \qquad\quad
  \le \epsilon t^{-\frac{K}{\theta}}
  +C_3t^{-\frac{K}{\theta}}\int_T^t E_{K,q}[F](s)\,ds
  \end{split}
  \end{equation}
  for all sufficiently large $t>T$ and all $0\le\ell\le K$. 
  In particular, if 
  $$
  \int_0^\infty E_{K,q}[F](s)\,ds<\infty,
  $$
  then 
  \begin{equation}
  \label{eq:3.13}
  \lim_{t\to\infty}
  t^{\frac{K}{\theta}}
  \biggr[t^{\frac{N}{\theta}(1-\frac{1}{q})}\|v(t)-w_2(t)\|_q+t^{-\frac{\ell}{\theta}}|||v(t)-w_2(t)|||_\ell\biggr]=0.
  \end{equation}
\end{itemize}
The constants $C_1$, $C_2$ and $C_3$ depend only on $N$, $\theta$ and $K$. 
\end{theorem}
We first prove assertion~(i) of Theorem~\ref{Theorem:3.1}.
\vspace{3pt}
\newline
{\bf Proof of Theorem~\ref{Theorem:3.1}~(i).}
Let $K\ge0$ and $\varphi\in L^\infty\cap L_K$. 
Put 
\begin{equation}
\label{eq:3.14}
\tilde{\varphi}(x):=\varphi(x)-\sum_{|\alpha|\le K}M_\alpha(\varphi,0)\psi_\alpha(x,0). 
\end{equation}
By \eqref{eq:3.7} we see that
\begin{equation*}
\int_{{\bf R}^N}x^\beta\tilde{\varphi}(x)\,dx=0\quad\mbox{for all $\beta\in{\bf M}_K$}.
\end{equation*}
This together with Theorem~\ref{Theorem:1.1}, \eqref{eq:3.3} and \eqref{eq:3.14} implies
\begin{eqnarray}
\label{eq:3.15}
 & & 
\sup_{t>0}\,t^{\frac{K}{\theta}}
\biggr[t^{\frac{N}{\theta}(1-\frac{1}{q})}\|S_\theta(t)\tilde{\varphi}\|_q
+t^{-\frac{\ell}{\theta}}|||S_\theta(t)\tilde{\varphi}|||_\ell\biggr]\\
\notag
 & & \qquad\quad
\le C|||\tilde{\varphi}|||_K\le C(\|\varphi\|_1+|||\varphi|||_K),\\
\label{eq:3.16}
 & & \lim_{t\to\infty}
t^{\frac{K}{\theta}}
\biggr[t^{\frac{N}{\theta}(1-\frac{1}{q})}\|S_\theta(t)\tilde{\varphi}\|_q
+t^{-\frac{\ell}{\theta}}|||S_\theta(t)\tilde{\varphi}|||_\ell\biggr]=0,
\end{eqnarray}
for $1\le q\le\infty$ and $0\le\ell\le K$.
On the other hand, 
it follows from \eqref{eq:3.4} and \eqref{eq:3.14} that
\begin{equation}
\label{eq:3.17}
[S_\theta(t)\tilde{\varphi}](x)
=[S_\theta(t)\varphi](x)-w_1(x,t),
\quad
x\in{\bf R}^N,\,\,t>0. 
\end{equation} 
By \eqref{eq:2.4}, \eqref{eq:3.1}, \eqref{eq:3.2}, \eqref{eq:3.3} and \eqref{eq:3.17} we obtain
\begin{equation}
\label{eq:3.18}
\begin{split}
\|S_\theta(t)\tilde{\varphi}\|_q & 
\le \|S_\theta(t)\varphi\|_q+\sum_{|\alpha|\le K}|M_\alpha(\varphi,0)|\|\Psi_{\alpha,\theta}(t:0)\|_q\\
 & \le C(\|\varphi\|_\infty+\|\varphi\|_1+|||\varphi|||_K),
\qquad t>0,
\end{split}
\end{equation}
for $1\le q\le\infty$.
This together with \eqref{eq:3.15} implies
\begin{equation}
\label{eq:3.19}
\begin{split}
|||S_\theta(t)\tilde{\varphi}|||_\ell
 & =(1+t)^{\frac{\ell}{\theta}}\int_{{\bf R}^N}\biggr(\frac{|x|}{(1+t)^{\frac{1}{\theta}}}\biggr)^\ell
\left|[S_\theta(t)\tilde{\varphi}](x)\right|\,dx\\
 & \le C(1+t)^{\frac{\ell}{\theta}}
\int_{{\bf R}^N}\biggr[1+\biggr(\frac{|x|}{(1+t)^{\frac{1}{\theta}}}\biggr)^K\biggr]
\left|[S_\theta(t)\tilde{\varphi}](x)\right|\,dx\\
 & =C(1+t)^{\frac{\ell}{\theta}}\|S_\theta(t)\tilde{\varphi}\|_1+C(1+t)^{\frac{\ell-K}{\theta}}|||S_\theta(t)\tilde{\varphi}|||_K\\
 & \le C(1+t)^{\frac{\ell}{\theta}}(\|\varphi\|_\infty+\|\varphi\|_1+|||\varphi|||_K),
 \qquad t>0,
\end{split}
\end{equation}
for $0\le\ell\le K$. 
Therefore, by \eqref{eq:3.15}, \eqref{eq:3.16}, \eqref{eq:3.17}, \eqref{eq:3.18} and \eqref{eq:3.19} 
we obtain \eqref{eq:3.8} and \eqref{eq:3.9}. 
Thus Theorem~\ref{Theorem:3.1}~(i) follows.
$\Box$
\vspace{5pt}

For the proof of Theorem~\ref{Theorem:3.1}~(ii), 
we prepare the following lemma, which is proved by a similar argument as in \cite[Lemma~2.2]{IKK02} 
with the aid of \eqref{eq:3.3}. 
\begin{lemma}
\label{Lemma:3.1} 
Assume the same assumptions as in Theorem~{\rm\ref{Theorem:3.1}~(ii)}. 
Let $1\le r\le q\le\infty$ and $0\le\ell\le K$. 
Set 
\begin{equation}
\label{eq:3.20}
\tilde{F}(x,t):=F(x,t)-\sum_{|\alpha|\le K}M_\alpha(F(t),t)\,\psi_\alpha(x,t). 
\end{equation}
Then there exists a constant $C$ such that 
\begin{equation*}
\begin{split}
 & \|\tilde{F}(t)\|_r\le C(1+t)^{-\frac{N}{\theta}(1-\frac{1}{r})-\frac{K}{\theta}}E_{K,q}[F](t),\\
 & |||\tilde{F}(t)|||_\ell\le C(1+t)^{-\frac{K-\ell}{\theta}}E_{K,q}[F](t),
\end{split} 
\end{equation*}
for $t>0$.
\end{lemma}
{\bf Proof of Theorem~\ref{Theorem:3.1}~(ii).}
Let $0\le\ell\le K$ and $1\le q\le\infty$. 
It follows from \eqref{eq:3.4} and \eqref{eq:3.20} that 
\begin{equation}
\label{eq:3.21}
\begin{split}
 & v(x,t)-w_2(x,t)
=\int_0^t[S_\theta(t-s)\tilde{F}(\cdot,s)](x)\,ds\\
 & \qquad\quad
 =\int_{t/2}^t[S_\theta(t-s)\tilde{F}(\cdot,s)](x)\,ds
 +\int_0^{t/2}[S_\theta(t-s)\tilde{F}(\cdot,s)](x)\,ds\\
 & \qquad\quad
 =: I_1(x,t)+I_2(x,t)
\end{split}
\end{equation}
for $(x,t)\in{\bf R}^N\times(0,\infty)$.
Furthermore, by \eqref{eq:3.7} we see that 
\begin{equation}
\label{eq:3.22}
\int_{{\bf R}^N}x^\beta\tilde{F}(x,s)\,dx=0\quad\mbox{for all $\beta\in{\bf M}_K$, $s>0$}. 
\end{equation} 
This together with Theorem~\ref{Theorem:1.1}~(a) with $K=\ell$, 
Lemma~\ref{Lemma:3.1} and \eqref{eq:2.4} implies that 
\begin{equation*}
\begin{split}
 & \|v(t)-w_2(t)\|_q+|||v(t)-w_2(t)|||_\ell\\
 & \le C\int_0^t\|\tilde{F}(s)\|_q\,ds
 +C\int_0^t|||\tilde{F}(s)|||_\ell\,ds
 \le C\int_0^t E_{K,q}[F](s)\,ds
\end{split}
\end{equation*}
for $t>0$. This implies \eqref{eq:3.11}. 
Similarly, we have
\begin{equation}
\label{eq:3.23}
\begin{split}
 & t^{\frac{N}{\theta}(1-\frac{1}{q})}\|I_1(t)\|_q+(1+t)^{-\frac{\ell}{\theta}}|||I_1(t)|||_\ell\\
 & \le Ct^{\frac{N}{\theta}(1-\frac{1}{q})}\int_{t/2}^t\|\tilde{F}(s)\|_q\,ds+C(1+t)^{-\frac{\ell}{\theta}}\int_{t/2}^t|||\tilde{F}(s)|||_\ell\,ds\\
 & \le Ct^{-\frac{K}{\theta}}\int_{t/2}^t |||\tilde{F}(s)|||_K\,ds
\le Ct^{-\frac{K}{\theta}}\int_{t/2}^t E_{K,q}[F](s)\,ds.
\end{split}
\end{equation}

Let $T>0$. 
By Theorem~\ref{Theorem:1.1}, Lemma~\ref{Lemma:3.1} and \eqref{eq:3.22} 
we apply the Lebesgue dominated convergence theorem to obtain
\begin{equation}
\label{eq:3.24}
\limsup_{t\to\infty}\,t^{\frac{K}{\theta}}
\int_0^T 
\left[t^{\frac{N}{\theta}(1-\frac{1}{q})}\|S_\theta(t-s)\tilde{F}(\cdot,s)\|_q
+t^{-\frac{\ell}{\theta}}|||S_\theta(t-s)\tilde{F}(\cdot,s)|||_\ell\right]\,ds=0.
\end{equation} 
Furthermore, 
by Theorem~\ref{Theorem:1.1}~(a), Lemma~\ref{Lemma:3.1} and \eqref{eq:3.22}
we obtain
\begin{equation}
\label{eq:3.25}
\begin{split}
 & t^{\frac{N}{\theta}(1-\frac{1}{q})}\int_T^{t/2}
 \|S_\theta(t-s)\tilde{F}(\cdot,s)\|_q\,ds
 +t^{-\frac{\ell}{\theta}}\int_T^{t/2}
 |||S_\theta(t-s)\tilde{F}(\cdot,s)|||_\ell\,ds\\
 & \le Ct^{\frac{N}{\theta}(1-\frac{1}{q})}\int_T^{t/2}
 (t-s)^{-\frac{N}{\theta}(1-\frac{1}{q})-\frac{K}{\theta}}
 |||\tilde{F}(s)|||_K\,ds\\
 & \qquad\qquad\qquad
 +Ct^{-\frac{\ell}{\theta}}\int_T^{t/2} 
(t-s)^{-\frac{K-\ell}{\theta}}|||\tilde{F}(s)|||_K\,ds\\
 & \le Ct^{-\frac{K}{\theta}}\int_T^{t/2} E_{K,q}[F(s)]\,ds
\end{split}
\end{equation} 
for $t\ge 2T$. 
Then, for any $\epsilon>0$, we deduce from \eqref{eq:3.24} and \eqref{eq:3.25} that
\begin{equation}
\label{eq:3.26}
t^{\frac{N}{\theta}(1-\frac{1}{q})}\|I_2(t)\|_q
+t^{-\frac{\ell}{\theta}}|||I_2(t)|||_\ell
\le\epsilon t^{-\frac{K}{\theta}}+Ct^{-\frac{K}{\theta}}\int_T^{t/2}E_{K,q}[F](s)\,ds
\end{equation} 
for all sufficiently large $t$.
Then, by \eqref{eq:3.21}, \eqref{eq:3.23} and \eqref{eq:3.26}
we obtain \eqref{eq:3.12}. 
Furthermore, 
\eqref{eq:3.13} immediately follows from \eqref{eq:3.12}. 
Therefore the proof of Theorem~\ref{Theorem:3.1}~(ii) is complete. 
Thus Theorem~\ref{Theorem:3.1} follows. 
$\Box$
\begin{Remark}
\label{Remark:3.1}
{\rm (i)} Let $0\le K<\theta$ and set $\psi(x)=G_\theta(x,1)$. 
It follows that $\psi(x,s)=G_\theta(x,1+s)$ and 
$$
\Psi_{\alpha,\theta}(x,t:s)=\frac{(-1)^{|\alpha|}}{\alpha!}\partial_x^\alpha[S_\theta(t)G_\theta(\cdot,1+s)](x)
=g_{\alpha,\theta}(x,1+s+t)
$$
for $x\in{\bf R}^N$, $t>0$ and $s\ge 0$. 
Furthermore, the functions $w_1$ and $w_2$ in Theorem~{\rm\ref{Theorem:3.1}} are represented by 
\begin{equation*}
\begin{split}
 & w_1(x,t)=\sum_{|\alpha|\le K}M_\alpha(\varphi,0)g_{\alpha,\theta}(x,t+1),\\
 & w_2(x,t)=\sum_{|\alpha|\le K}
\biggr[\int_0^t M_\alpha(F(s),s)\,ds\biggr]g_{\alpha,\theta}(x,t+1),
\end{split}
\end{equation*}
respectively. Then, by Theorem~{\rm\ref{Theorem:3.1}} 
we obtain the same conclusions as in Theorems~{\rm 1.1} and {\rm 1.2} in {\rm\cite{IKK02}}. 
\vspace{3pt}
\newline
{\rm (ii)} 
In the case $K\ge\theta$, since $G_\theta(\cdot,1)\not\in L_K$, 
the argument in Remark~{\rm \ref{Remark:3.1}}~{\rm (i)} is not applicable. 
\vspace{3pt}
\newline
{\rm (iii)} 
It follows from \eqref{eq:3.1} and \eqref{eq:3.2} that $\psi_\alpha\in L^\infty\cap L_K$.  
This fact enables us to obtain a better asymptotic expansion $w_1$ of $S_\theta(t)\varphi$ than 
$$
\sum_{|\alpha|\le K}
\biggr[\int_{{\bf R}^N}x^\alpha \varphi(x)\,dx\biggr]g_{\alpha,\theta}(x,t)
$$
{\rm(}compare \eqref{eq:3.8} with \eqref{eq:1.10}{\rm)} and to define the function $w_2$ in 
assertion~{\rm (ii)} of Theorem~{\rm\ref{Theorem:3.1}}. 
\end{Remark}
%
%%%%%%%%%%%%%%%%%%%%%%%%%%%%%%%%%%%%%
%%%%%%%%%%%%%%%%%%%%%%%%%%%%%%%%%%%%%
\section{Nonlinear fractional diffusion equation}
%%%%%%%%%%%%%%%%%%%%%%%%%%%%%%%%%%%%%
%%%%%%%%%%%%%%%%%%%%%%%%%%%%%%%%%%%%%
Let $K\ge 0$. Assume that a radially symmetric smooth function $\psi$ satisfies \eqref{eq:3.1} and 
\begin{equation}
\label{eq:4.1}
|(\partial_x^\alpha\psi)(x)|\le CG_\theta(x,1),\qquad \alpha\in{\bf M}_K,
\end{equation}
for $x\in{\bf R}^N$. 
We apply the arguments in the previous sections 
to study the asymptotic expansions of the solution~$u$ of \eqref{eq:1.1} satisfying \eqref{eq:1.4}. 
Let
$$
\mathcal{L}_K^*:=L_K+
\left\{\phi\in L^\infty\,:\,
\underset{x\in{\bf R}^N}{\mbox{ess sup}}
\,(1+|x|)^{N+\theta}|\phi(x)|<\infty\right\}.
$$
\begin{theorem}
\label{Theorem:4.1}
Let $N\ge 1$, $0<\theta<2$, $p>1+\theta/N$ and $K\ge 0$. Assume $\varphi\in L^\infty\cap \mathcal{L}_K^*$. 
Let $u$ be a global in time solution of \eqref{eq:1.1} such that 
\begin{equation}
\label{eq:4.2}
\sup_{t>0}\,(1+t)^{\frac{N}{\theta}}\|u(t)\|_\infty<\infty. 
\end{equation}
Futhermore, let $\psi$ be a radially symmetric smooth function in ${\bf R}^N$ satisfying \eqref{eq:3.1} and \eqref{eq:4.1}. 
Define
$$
U_0(x,t):=
[S_\theta(t)\varphi](x)+\sum_{|\alpha|\le K}
\int_0^tM_\alpha(F(s),s)\Psi_{\alpha,\theta}(x,t-s:s)\,ds,
$$
where $M_\alpha(F(s),s)$ and $\Psi_{\alpha,\theta}(x,t-s:s)$ are as in Section~{\rm 3}.
Assume $K+N<p(N+\theta)$. 
Then, for any $1\le q\le\infty$ and $0\le\ell\le K$, 
\begin{equation}
\label{eq:4.3}
\sup_{t>0}\,
\left[(1+t)^{\frac{N}{\theta}(1-\frac{1}{q})}\left\|u(t)-U_0(t)\right\|_q
 +(1+t)^{-\frac{\ell}{\theta}}
\left|\left|\left|u(t)-U_0(t)\right|\right|\right|_\ell\right]<\infty
\end{equation}
and
\begin{align}
 & t^{\frac{N}{\theta}(1-\frac{1}{q})}\left\|u(t)-U_0(t)\right\|_q
 +t^{-\frac{\ell}{\theta}}
\left|\left|\left|u(t)-U_0(t)\right|\right|\right|_\ell\vspace{3pt}\notag\\
 & \qquad\qquad\qquad
 =\left\{
\begin{array}{ll}
o(t^{-\frac{K}{\theta}})+O(t^{-A_p+1})  & \mbox{if}\quad A_p-1\not=K/\theta,\vspace{3pt}\\
O(t^{-\frac{K}{\theta}}\log t)  & \mbox{if}\quad A_p-1=K/\theta,\\
\end{array}
\right.
\label{eq:4.4}
\end{align}
as $t\to\infty$, where $A_p:=N(p-1)/\theta>1$. 
\end{theorem}
For the proof of Theorem~\ref{Theorem:4.1}, 
we prepare the following two lemmas on ${\mathcal L}_K^*$. 
\begin{lemma}
\label{Lemma:4.1}
Let $0<\theta<2$ and $K\ge0$. Then 
\begin{itemize}
  \item[{\rm (a)}]
  ${\mathcal L}_{K,0}\subset{\mathcal L}_K^*$;
  \item[{\rm (b)}]
  $\displaystyle{\bigcup_{t>0}\,
  \left\{S_\theta(t)\varphi\,:\,\varphi\in{\mathcal L}_K^*\right\}\subset{\mathcal L}_K^*}$;
  \item[{\rm (c)}]
  Let $\varphi\in{\mathcal L}_K^*$ and let $\phi$ be a measurable function in ${\bf R}^N$ 
  such that $|\phi(x)|\le|\varphi(x)|$ 
  for almost all $x\in{\bf R}^N$. Then $\phi\in{\mathcal L}_K^*$. 
\end{itemize}
\end{lemma}
{\bf Proof.}
Assertion~(a) follows from Lemma~\ref{Lemma:2.1}. 
Let $\varphi\in{\mathcal L}_K^*$. Then we can find 
$\varphi_1\in L_K$ and $\varphi_2\in L^\infty$ with 
\begin{equation}
\label{eq:4.5}
\underset{x\in{\bf R}^N}{\mbox{ess sup}}
\,(1+|x|)^{N+\theta}|\varphi_2(x)|<\infty
\end{equation}
such that $\varphi=\varphi_1+\varphi_2$. 

We prove assertion~(b). 
It follows from Theorem~\ref{Theorem:1.2}~(b) that 
\begin{equation}
\label{eq:4.6}
\{S_\theta(t)\varphi_1\,:\,t>0\}\subset{\mathcal L}_{K,0}. 
\end{equation}
By ({\bf G})-(iii) and \eqref{eq:4.5} we have 
$$
|\varphi_2(x)|\le CG_\theta(x,1)
$$
for almost all $x\in{\bf R}^N$.
Then 
$$
|[S_\theta(t)\varphi_2](x)|\le [S_\theta(t)|\varphi_2|](x)
\le CG_\theta(x,1+t),\quad x\in{\bf R}^N,\,\,\,t>0,
$$
which together with \eqref{eq:2.2} implies that 
\begin{equation}
\label{eq:4.7}
\sup_{x\in{\bf R}^N}\,(1+|x|)^{N+\theta}
|[S_\theta(t)\varphi_2](x)|<\infty. 
\end{equation}
Therefore, by assertion~(a), \eqref{eq:4.6} and \eqref{eq:4.7} we obtain 
$\left\{S_\theta(t)\varphi\,:\,\varphi\in{\mathcal L}_K^*\right\}\subset{\mathcal L}_K^*$ for any $t>0$. 
Thus assertion~(b) follows. 

We prove assertion~(c). 
Set $$
\phi^\pm:=\max\{\pm\phi,0\},
\qquad
E^\pm:=\{x\in{\bf R}^N\,:\,\phi^\pm(x)>|\varphi_1(x)|\}.
$$ 
Since $|\phi(x)|\le|\varphi(x)|\le|\varphi_1(x)|+|\varphi_2(x)|$ for almost all $x\in{\bf R}^N$, 
we see that  
\begin{equation}
\label{eq:4.8}
0\le (\phi^\pm(x)-|\varphi_1(x)|)\chi_{E^\pm}(x)\le|\varphi_2(x)|\le C(1+|x|)^{-N-\theta}
\end{equation}
for almost all $x\in{\bf R}^N$, 
where $\chi_{E^\pm}$ are the characteristic functions of $E^\pm$, respectively. 
Furthermore, 
since $|\varphi_1|\chi_{E^\pm}\le|\varphi_1|$ and $\phi^\pm\chi_{{\bf R}^N\setminus E^\pm}\le|\varphi_1|$,
we have 
\begin{equation}
\label{eq:4.9}
|\varphi_1|\chi_{E^\pm}+\phi^\pm\chi_{{\bf R}^N\setminus E^\pm}\in L_K. 
\end{equation}
Therefore, by \eqref{eq:4.8} and \eqref{eq:4.9} 
we obtain 
$$
\phi^\pm=(\phi^\pm-|\varphi_1|)\chi_{E^\pm}+|\varphi_1|\chi_{E^\pm}+\phi^\pm\chi_{{\bf R}^N\setminus E^\pm}
\in{\mathcal L}_K^*.
$$
This implies that $\phi\in{\mathcal L}_K^*$ and assertion~(c) follows. 
Therefore the proof of Lemma~\ref{Lemma:4.1} is complete. 
$\Box$
\begin{lemma}
\label{Lemma:4.2}
Let $u$ be a global in time solution of \eqref{eq:1.1} with $\varphi\in L^\infty\cap \mathcal{L}_K^*$ 
and $p>1+\theta/N$. 
Assume \eqref{eq:4.2}. 
Then 
\begin{equation}
\label{eq:4.10}
u(t)\in L^\infty\cap\mathcal{L}_K^*,\qquad t>0. 
\end{equation}
Furthermore, if $K+N<p(N+\theta)$,
then there exists a constant $C$ such that 
\begin{equation}
\label{eq:4.11}
|||u(t)^p|||_K\le C(1+t)^{-A_p+\frac{K}{\theta}},\qquad t>0.
\end{equation}
\end{lemma}
{\bf Proof.}
We prove \eqref{eq:4.10}.
It follows from \eqref{eq:4.2} that 
$$
\partial_t u+(-\Delta)^{\frac{\theta}{2}}u\le C(1+t)^{-A_p}u\quad\mbox{in}\quad{\bf R}^N\times(0,\infty).
$$
Since $A_p>1$, the comparison principle implies that 
\begin{equation}
\label{eq:4.12}
|u(x,t)|\le C[S_\theta(t)|\varphi|](x),
\qquad x\in{\bf R}^N,\,\,\,t>0.
\end{equation}
It follows from $\varphi\in\mathcal{L}_K^*$ that $|\varphi|\in{\mathcal L}_K^*$.
This together with Lemma~\ref{Lemma:4.1}~(b) yields $S_\theta(t)|\varphi|\in{\mathcal L}_K^*$.
Then, by Lemma~\ref{Lemma:4.1}~(c) and \eqref{eq:4.12} 
we see that $u(t)\in{\mathcal L}_K^*$ for $t>0$.
This together with \eqref{eq:4.2} implies \eqref{eq:4.10}. 

We prove \eqref{eq:4.11}.
We assume that $K+N<p(N+\theta)$.
By \eqref{eq:2.2} we see that 
$G_\theta(\cdot,t)^p\in L_K$ and 
\begin{equation}
\label{eq:4.13}
|||G_\theta(t)^p|||_K\le Ct^{-A_p+\frac{K}{\theta}},\qquad t>0. 
\end{equation}
Let $\varphi\in L^\infty\cap{\mathcal L}_K^*$. 
Then we can find 
$\varphi_1\in L^\infty\cap L_K$ and $\varphi_2\in L^\infty$ with 
\begin{equation}
\label{eq:4.14}
|\varphi_2(x)|\le CG_\theta(x,1)\qquad \mbox{for almost all}\,\,\,x\in{\bf R}^N,
\end{equation}
such that $\varphi=\varphi_1+\varphi_2$.
By  $\varphi_1\in L_K$ we define
\begin{equation}
\label{eq:4.15}
v_1(x,t):=[S_\theta(t)|\varphi_1|](x)
-\sum_{|\alpha|\le K}M_\alpha(|\varphi_1|,0)\Psi_{\alpha,\theta}(x,t:0),
\end{equation}
where $\Psi_{\alpha,\theta}$ is the function given by \eqref{eq:3.2}.
Then Theorem~\ref{Theorem:3.1}~(i) implies that
\begin{equation}
\label{eq:4.16}
|||v_1(t)^p|||_K\le \|v_1(t)\|_\infty^{p-1}|||v_1(t)|||_K
\le C(1+t)^{-A_p-\frac{K(p-1)}{\theta}},
\quad t>0.
\end{equation}
On the other hand, 
by \eqref{eq:3.4} and \eqref{eq:4.1} we have
\begin{equation*}
\begin{split}
\Psi_{\alpha,\theta}(x,t:0)
 & =[S_\theta(t)\psi_\alpha(0)](x)
\le[S_\theta(t)|\psi_\alpha(0)|](x)\\
 & \le C[S_\theta(t)G_\theta(\cdot,1)](x)
\le CG_\theta(x,t+1)
\end{split}
\end{equation*}
for $x\in{\bf R}^N$, $t>0$ and $\alpha\in{\bf M}_K$.
This together with \eqref{eq:4.13}, \eqref{eq:4.15} and \eqref{eq:4.16}
yields 
\begin{equation}
\label{eq:4.17}
\begin{split}
|||(S_\theta(t)|\varphi_1|)^p|||_K
&
\le C|||v_1(t)^p|||_K+C\sum_{|\alpha|\le K}|||\Psi_{\alpha,\theta}(t:0)^p|||_K
\\
&
\le C(1+t)^{-A_p-\frac{K(p-1)}{\theta}}+C|||G_\theta(t+1)^p|||_K
\le C(1+t)^{-A_p+\frac{K}{\theta}}
\end{split}
\end{equation}
for $t>0$.
Furthermore, by \eqref{eq:4.13} and \eqref{eq:4.14} we obtain
\begin{equation}
\label{eq:4.18}
 |||(S_\theta(t)|\varphi_2|)^p|||_K\le C|||[S_\theta(t) G_\theta(1)]^p|||_K
 =C|||G_\theta(t+1)^p|||_K
 \le C(1+t)^{-A_p+\frac{K}{\theta}}
\end{equation}
for $t>0$. 
On the other hand, by \eqref{eq:4.12} we see that 
$$
|u(x,t)|^p\le C[S_\theta(t)|\varphi|](x)^p
\le C[S_\theta(t)|\varphi_1|](x)^p+C[S_\theta(t)|\varphi_2|](x)^p
$$
for $x\in{\bf R}^N$ and $t>0$. 
This together with \eqref{eq:4.17} and \eqref{eq:4.18} implies that 
$$
|||u(t)^p|||_K\le C(1+t)^{-A_p+\frac{K}{\theta}},
\qquad t>0. 
$$
Therefore we obtain \eqref{eq:4.11}, 
and the proof of Lemma~\ref{Lemma:4.2} is complete.
$\Box$
\vspace{5pt}

Now we are ready to prove Theorem~\ref{Theorem:4.1}.
\vspace{5pt}
\newline
{\bf Proof of Theorem~\ref{Theorem:4.1}.}
Let $1\le q\le\infty$ and $0\le\ell\le K$. 
Since $\varphi\in{\mathcal L}_K^*$, 
by \eqref{eq:2.4}, \eqref{eq:4.2} and \eqref{eq:4.12} we have
\begin{equation}
\label{eq:4.19}
\|F(u(t))\|_q\le 
\|u(t)\|_\infty^{p-1}\|u(t)\|_q\le C (1+t)^{-A_p-\frac{N}{\theta}(1-\frac{1}{q})},
\qquad t>0.
\end{equation}
Then, by \eqref{eq:3.10}, \eqref{eq:4.11} and \eqref{eq:4.19} we have
\begin{equation}
\label{eq:4.20}
E_{K,q}[F](t)\le C(1+t)^{-A_p+\frac{K}{\theta}},\qquad t>0.
\end{equation}
On the other hand, it follows from Theorem~\ref{Theorem:3.1}~(ii) that 
$$
\|u(t)-U_0(t)\|_q+|||u(t)-U_0(t)|||_\ell\le C\int_0^t E_{K,q}[F](s)\,ds,
\qquad t>0. 
$$
Furthermore, 
for any $T>0$ and $\epsilon>0$, we obtain
$$
t^{\frac{N}{\theta}(1-\frac{1}{q})}\|u(t)-U_0(t)\|_q+t^{-\frac{\ell}{\theta}}|||u(t)-U_0(t)|||_\ell
\le \epsilon t^{-\frac{K}{\theta}}+
C_*t^{-\frac{K}{\theta}}\int_T^t E_{K,q}[F](s)\,ds
$$
for all sufficiently large $t>T$, where $C_*$ is a constant independent of $T$ and $\epsilon$. 
These together with \eqref{eq:4.20} imply \eqref{eq:4.3} and \eqref{eq:4.4}. 
Thus Theorem~\ref{Theorem:4.1} follows.
$\Box$
\vspace{7pt}

Finally, combining the arguments in \cite[Theorem~3.1]{IK02}, 
we obtain the following theorem.  
\begin{theorem}
\label{Theorem:4.2}
Assume the same conditions as in Theorem~{\rm\ref{Theorem:4.1}}. 
Define $U_n=U_n(x,t)$ $(n=1,2,\dots)$ inductively by
\begin{equation}
\label{eq:4.21}
\begin{split}
U_n(x,t) & := 
U_0(x,t)+\int_0^t [S_\theta(t-s)F_{n-1}(\cdot,s)](x)\,ds\\
 & \quad
-\sum_{|\alpha|\le K}\int_0^tM_\alpha(F_{n-1}(s),s)\Psi_{\alpha,\theta}(x,t-s:s)\,ds,
\end{split}
\end{equation}
where $n=1,2,\dots$ and $F_{n-1}(x,t):=F(U_{n-1}(x,t))$. 
Then, 
for any $1\le q\le\infty$, $0\le \ell\le K$, and $n=0,1,\dots$,
\begin{equation}
\label{eq:4.22}
\sup_{t>0}\left\{(1+t)^{\frac{N}{\theta}(1-\frac{1}{q})}\|u(t)-U_n(t)\|_q
+(1+t)^{-\frac{K}{\theta}}|||u(t)-U_n(t)|||_K\right\}<\infty
\end{equation}
and
\begin{align}
 & t^{\frac{N}{\theta}(1-\frac{1}{q})}\left\|u(t)-U_n(t)\right\|_q
 +t^{-\frac{\ell}{\theta}}\left|\left|\left|u(t)-U_n(t)\right|\right|\right|_\ell\vspace{3pt}\notag\\
 & \qquad
 =\left\{
\begin{array}{ll}
o(t^{-\frac{K}{\theta}})+O(t^{-(n+1)(A_p-1)})  & \mbox{if}\quad (n+1)(A_p-1)\not=K/\theta,\vspace{3pt}\\
O(t^{-\frac{K}{\theta}}\log t)  & \mbox{if}\quad (n+1)(A_p-1)=K/\theta,\\
\end{array}
\right.
\label{eq:4.23}
\end{align}
as $t\to\infty$. 
\end{theorem}
{\bf Proof of Theorem~\ref{Theorem:4.2}.}
By Theorem~\ref{Theorem:4.1} we have \eqref{eq:4.22} and \eqref{eq:4.23} with $n=0$.
Assume that \eqref{eq:4.22} and \eqref{eq:4.23} hold for some $n=m\in\{0,1,2,\dots\}$. 
Set 
$$
\hat{F}(x,t):=F(u(x,t))-F(U_m(x,t)).
$$
Then it follows from \eqref{eq:4.21} that 
\begin{equation}
\label{eq:4.24}
u(x,t)-U_{m+1}(x,t)
 =\hat{I}(x,t)
\end{equation}
for $x\in{\bf R}^N$ and $t>0$,
where
$$
\hat{I}(x,t):=
\int_0^t[S_\theta(t-s)\hat{F}(\cdot,s)](x)\,ds
 -\sum_{|\alpha|\le K}\int_0^tM_\alpha(\hat{F}(s),s)\Psi_{\alpha,\theta}(x,t-s:s)\,ds.
$$
By \eqref{eq:4.2} and \eqref{eq:4.22} 
we apply the mean value theorem to obtain 
$$
E_{K,q}[\hat{F}](t)\le C(1+t)^{-A_p}E_{K,q}[u-U_m](t),\qquad t>0,
$$
where $1\le q\le\infty$. 
This together with \eqref{eq:3.10}, \eqref{eq:4.22} and \eqref{eq:4.23} with $n=m$ implies that 
$$
  E_{K,q}[\hat{F}](t)\le C(1+t)^{-A_p+\frac{K}{\theta}},\qquad t>0,$$
and 
$$
 E_{K,q}[\hat{F}](t)=\left\{
\begin{array}{l}
o(t^{-A_p})+O(t^{-A_p-(m+1)(A_p-1)+\frac{K}{\theta}})\vspace{3pt}\\
\qquad\qquad\qquad
\mbox{if}\quad (m+1)(A_p-1)\not=K/\theta,\vspace{3pt}\\
O(t^{-A_p}\log t)\vspace{3pt}\\
\qquad\qquad\qquad
\mbox{if}\quad (m+1)(A_p-1)=K/\theta,
\end{array}
\right.
$$
as $t\to\infty$. Then we deduce from Theorem~\ref{Theorem:3.1}~(ii) that 
\begin{equation}
\label{eq:4.25}
\|\hat{I}(t)\|_q+|||\hat{I}(t)|||_\ell\le C\int_0^t(1+s)^{-A_p+\frac{K}{\theta}}\,ds,\qquad t>0
\end{equation}
and
\begin{equation}
\label{eq:4.26}
\begin{split}
 & t^{\frac{N}{\theta}(1-\frac{1}{q})}\|\hat{I}(t)\|_q+t^{-\frac{\ell}{\theta}}|||\hat{I}(t)|||_\ell\\
 & =\left\{
\begin{array}{ll}
o(t^{-\frac{K}{\theta}})+O(t^{-(m+2)(A_p-1)})  & \mbox{if}\quad (m+2)(A_p-1)\not=K/\theta,\vspace{3pt}\\
O(t^{-\frac{K}{\theta}}\log t)  & \mbox{if}\quad (m+2)(A_p-1)=K/\theta,\\
\end{array}
\right.
\end{split}
\end{equation}
as $t\to\infty$, for any $1\le q\le\infty$ and $0\le\ell\le K$. 
Therefore, by \eqref{eq:4.24}, \eqref{eq:4.25} and \eqref{eq:4.26}
we obtain \eqref{eq:4.22} and \eqref{eq:4.23} with $n=m+1$.
This means that \eqref{eq:4.22} and \eqref{eq:4.23} hold for all $n=0,1,\dots$. 
Thus Theorem~\ref{Theorem:4.2} follows. 
$\Box$
%%%%%%%%%%%%%%%%%%%%%%%%%%%%%%%%%%%%%%
%%%%%%%%%%%%    references    %%%%%%%%%%%%%%%%%%
%%%%%%%%%%%%%%%%%%%%%%%%%%%%%%%%%%%%%%
\bibliographystyle{amsplain}

\begin{thebibliography}{10}

\bibitem{AF}
H. Amann and M. Fila, 
A Fujita-type theorem for the Laplace equation with a dynamical boundary condition,
Acta Math. Univ. Comenianae {\bf 66} (1997), 321--328.
     
\bibitem{BT}  
K. Bogdan and T. Jakubowski, 
Estimates of heat kernel of fractional Laplacian perturbed by gradient operators, 
Comm. Math. Phys.  {\bf 271}  (2007), 179--198. 

\bibitem{BK} 
L. Brandolese and G. Karch, 
Far field asymptotics of solutions to convection equation with anomalous diffusion, 
J. Evol. Equ. {\bf 8} (2008), 307--326.

\bibitem{CK}
I. Chavel and L. Karp, 
Movement of hot spots in Riemannian manifolds,
J. Analyse Math. {\bf 55} (1990), 271--286. 

\bibitem{DK} 
J. Dolbeault and G. Karch, 
Large time behavior of solutions to nonhomogeneous diffusion equations,  
Banach Center Publ. {\bf 74} (2006), 113--147.

\bibitem{FIK} 
M. Fila, K. Ishige and T. Kawakami, 
Convergence to the Poisson kernel for the Laplace equation with a nonlinear dynamical boundary condition,
Commun. Pure Appl. Anal. {\bf 11} (2012), 1285--1301.
	
\bibitem{FK} 
A. Fino and G. Karch,
Decay of mass for nonlinear equation with fractional Laplacian, 
Monatsh. Math. {\bf 160} (2010), 375--384.

\bibitem{FM} 
Y. Fujigaki and T. Miyakawa, 
Asymptotic profiles of nonstationary incompressible Navier-Stokes flows in the whole space, 
SIAM J. Math. Anal. {\bf 33} (2001), 523--544.
 
\bibitem{FI}  
Y. Fujishima and K. Ishige, 
Blow-up for a semilinear parabolic equation with large diffusion on ${\bf R}^N$, 
J. Differential Equations {\bf 250} (2011), 2508--2543. 

%\bibitem{GV}
%A. Gmira and L. V\'eron, 
%Large time behaviour of the solutions of a semilinear parabolic equation in ${\bf R}^N$, 
%J. Differential Equations {\bf 53} (1984), 258--276. 

\bibitem{HKN} 
N. Hayashi, E. I. Kaikina and P. I. Naumkin,
Asymptotics for fractional nonlinear heat equations, 
J. London Math. Soc. {\bf 72} (2005), 663--688. 

\bibitem{HI}
K. Hisa and K. Ishige, 
Existence of solutions for a fractional semilinear parabolic equation with singular initial data, 
preprint. 

\bibitem{IIK}
K. Ishige, M. Ishiwata and T. Kawakami, 
The decay of the solutions for the heat equation with a potential, 
Indiana Univ. Math. J. {\bf 58} (2009), 2673--2708. 

\bibitem{IKabe} 
K. Ishige and Y. Kabeya, 
$L^p$ norms of nonnegative Schr\"odinger heat semigroup and the large time behavior of hot spots, 
J. Funct. Anal. {\bf 262} (2012), 2695--2733.

\bibitem{IK01}
K. Ishige and T. Kawakami,
Asymptotic behavior of solutions for some semilinear heat equations in ${\bf R}^N$, 
Commun. Pure Appl. Anal. {\bf 8} (2009), 1351--1371.

\bibitem{IK02}
K. Ishige and T. Kawakami, 
Refined asymptotic profiles for a semilinear heat equation, 
Math. Ann. {\bf 353} (2012), 161--192.
	
\bibitem{IK03}
K. Ishige and T. Kawakami, 
Asymptotic expansions of solutions of the Cauchy problem for nonlinear parabolic equations, 
 J. Anal. Math. {\bf 121} (2013), 317--351. 

\bibitem{IKK01}
K. Ishige, K. Kawakami and K. Kobayashi,
Global solutions for a nonlinear integral equation with a generalized heat kernel, 
Discrete Contin. Dyn. Syst. Ser. S. {\bf 7} (2014), 767--783.

\bibitem{IKK02}
K. Ishige, T. Kawakami and K. Kobayashi, 
Asymptotics for a nonlinear integral equation with a generalized heat kernel, 
J. Evol. Equ. {\bf 14} (2014), 749--777. 
 
\bibitem{IKoba}
K. Ishige and K. Kobayashi, 
Convection-diffusion equation with absorption and non-decaying initial data,
J. Differential Equations {\bf 254} (2013), 1247--1268. 
	
\bibitem{Iwa}
T. Iwabuchi, 
Global solutions for the critical Burgers equation in the Besov spaces and the large time behavior,
Ann. Inst. H. Poincar\'e Anal. Non Lineair\'e {\bf 32} (2015), 687--713. 

\bibitem{KP} 
S. Kamin and L. A. Peletier, 
Large time behaviour of solutions of the heat equation with absorption, 
Ann. Scuola Norm. Sup. Pisa Cl. Sci. (4) {\bf 12} (1985), 393--408.

\bibitem{OY} 
T. Ogawa and M. Yamamoto, 
Asymptotic behavior of solutions to drift-diffusion system with generalized dissipation, 
Math. Models Methods Appl. Sci. {\bf 19} (2009), 939--967.

\bibitem{SS} 
Y. Shibata and S. Shimizu, 
A decay property of the Fourier transform and its application to the Stokes problem, 
J. Math. Fluid Mech. {\bf 3} (2001), 213--230.

\bibitem{S} 
S. Sugitani,
On nonexistence of global solutions for some nonlinear integral equations,
Osaka J. Math. {\bf 12} (1975), 45--51.
     
\bibitem{T} 
J. Taskinen, 
Asymptotical behaviour of a class of semilinear diffusion equations, 
J. Evol. Equ. {\bf 7} (2007), 429--447.

\bibitem{Y03} 
M. Yamamoto, 
Asymptotic expansion of solutions to the nonlinear 
dissipative equation with the anomalous diffusion, 
J. Math. Anal. Appl. {\bf 427}  (2015), 1027--1069.

\bibitem{YS} 
M. Yamamoto and Y. Sugiyama, 
Asymptotic expansion of solutions to the drift-diffusion equation with fractional dissipation, 
Nonlinear Anal. {\bf 141} (2016), 57--87.

\end{thebibliography}

%%%%%%%%%%%%%%%%%%%%%%%%%%%%%%%%%%%%%%%%%
\end{document}